\documentclass[12pt]{article}
\usepackage{bm, amsmath, amssymb, amsthm}

\newcommand{\R}{{\mathbb R}}
\newcommand{\Ric}{\mbox{Ric \,}}
\newcommand{\hess}{\mbox{Hess\,}}
\newcommand{\Sf}{{\mathbb S}}

\newtheorem{theorem}{Theorem}
\newtheorem{proposition}[theorem]{Proposition}
\newtheorem{lemma}[theorem]{Lemma}

\newtheorem{corollary}[theorem]{Corollary}

\newtheorem{remark}[theorem]{Remark}

\newcommand{\spa}{\mbox{span}}
\newcommand{\Q}{\mathbb{Q}}
\def\<{\langle}
\def\>{\rangle}
\def\d{\partial}
\def\be{\begin{equation} }
\def\ee{\end{equation} }
\newcommand\bea{\begin{eqnarray*}}
\newcommand\eea{\end{eqnarray*}}
\begin{document}

\title{Hypersurfaces of two space forms and conformally flat hypersurfaces}

\author{S. Canevari and R. Tojeiro}
\date{}
\maketitle

\begin{abstract}
We address the problem of determining the  hypersurfaces $f\colon M^{n} \to \Q_s^{n+1}(c)$ with dimension $n\geq 3$ of a 
pseudo-Riemannian space form of dimension $n+1$, constant curvature $c$ 
and  index $s\in \{0, 1\}$ for which  there exists another isometric immersion  $\tilde{f}\colon M^{n} \to \Q^{n+1}_{\tilde s}(\tilde{c})$  
with 
$\tilde{c}\neq c$. For $n\geq 4$, we provide a complete solution by extending results for  $s=0=\tilde s$ by do Carmo and Dajczer \cite{dcd} and by Dajczer and the second author \cite{dt1}.  
Our main results are for the most interesting case $n=3$, and these are new even in the Riemannian case $s=0=\tilde s$. 
In particular, we characterize the  solutions  that have dimension $n=3$ and three distinct
principal curvatures.  We show that these are closely related to  conformally flat hypersurfaces of $\Q_s^{4}(c)$ with three distinct principal curvatures, and we obtain a similar characterization of the latter that improves a theorem by Hertrich-Jeromin \cite{h-j}. 
We also derive a Ribaucour transformation for both classes of hypersurfaces,
which gives a process to produce a  family of new elements of those classes, starting from a given one, in terms of solutions of a linear system of PDE's. This enables us to construct explicit examples of three-dimensional solutions of the problem, as well as new explicit examples of three-dimensional conformally flat hypersurfaces that have three distinct principal curvatures.
\end{abstract}



We denote by $\Q_s^N(c)$   a pseudo-Riemannian space form of dimension $N$, constant sectional curvature $c$ 
and  index $s\in \{0, 1\}$, that is, $\Q_s^N(c)$ is either a Riemannian or Lorentzian space-form of constant curvature $c$, corresponding to  $s=0$ or $s=1$, respectively.
By a hypersurface $f\colon M^{n} \to \Q_s^{n+1}(c)$ we always mean  an isometric immersion of a \emph{Riemannian} manifold
$M^n$ of dimension $n$ into $\Q_s^{n+1}(c)$, thus $f$ is a \emph{space-like} hypersurface if $s=1$. 

One of the main purposes of  this paper is to address the following\vspace{1ex}\\
\emph{ Problem $*$: For which  hypersurfaces $f\colon M^{n} \to \Q_s^{n+1}(c)$ of dimension $n\geq 3$ does there exist another isometric immersion  $\tilde{f}\colon M^{n} \to \Q^{n+1}_{\tilde s}(\tilde{c})$  
with $\tilde{c}\neq c$?}\vspace{1ex}\\ 
This problem was studied for $s=0=\tilde s$ and $n\geq 4$ by do Carmo and Dajczer in \cite{dcd}, and by Dajczer and the second author in \cite{dt1}.
Some partial results in the most interesting case $n=3$  were also obtained in \cite{dt1}.
Including Lorentzian ambient space forms in our study of Problem $*$ was motivated by our investigation in \cite{ct} of  submanifolds of codimension two and constant curvature  
$c\in (0,1)$ of  $\Sf^{5}\times \mathbb{R}$,  which turned out to be related to  hypersurfaces $f\colon M^{3} \to \Sf^4$ 
 for which $M^3$ also admits an isometric immersion 
into the Lorentz space $\mathbb{R}_1^4=\Q_1^{4}(0)$. \vspace{1ex}


We first state our results for  the case $n\geq 4$. The next one extends a theorem due to  do Carmo and Dajczer 
\cite{dcd} in the case  $s=0=\tilde s$. 
Here and in the sequel, for $s, \tilde s\in \{0, 1\}$ we denote $\epsilon=-2s+1$ and $\tilde \epsilon=-2\tilde s+1$.

\begin{theorem}\label{ngeq4}
Let $f\colon M^{n} \to \Q_s^{n+1}(c)$ be a hypersurface of dimension $n\geq 4$. If there exists another isometric immersion  $\tilde{f}\colon M^{n} \to \Q^{n+1}_{\tilde s}(\tilde{c})$  
with $\tilde c \neq {c}$,  
then  $c< \tilde{c}$  if  $s=0$ and $\tilde s=1$ (respectively, $c > \tilde{c}$ if 
  $s=1$ and $\tilde s=0$) and  $f$  has a principal curvature $\lambda$ of multiplicity at least $n-1$ everywhere. Moreover,
  at any $x\in M^n$ the following holds:
  \begin{itemize}
  \item[$(i)$] if $\lambda=0$  or $f$ is umbilical  with $c+\epsilon\lambda^2\neq \tilde c$, then $\tilde f$ is umbilical;
  \item[$(ii)$] if $f$ is umbilical  and $c+\epsilon\lambda^2=\tilde c$, then $0$ is a principal curvature of $\tilde f$ 
  with  multiplicity at least $n-1$;
  \item[$(iii)$]  if $\lambda\neq 0$  with multiplicity $n-1$,  then $\tilde f$ has also a principal curvature $\tilde \lambda$  with the same eigenspace as $\lambda$.
  \end{itemize}
\end{theorem}

Thus, Problem $*$ has no solutions if $n\geq 4$ and either $c>\tilde c$, $s=0$ and $\tilde s=1$ or $c<\tilde c$, $s=1$ and $\tilde s=0$,
while, in the remaining cases, having a principal curvature of multiplicity  at least $n-1$  is a necessary condition for a solution.
In those cases, having a principal curvature of \emph{constant} multiplicity $n$ or $n-1$  is also sufficient for simply connected hypersurfaces.

\begin{theorem}\label{conv}
Let $f\colon M^{n} \to \Q_s^{n+1}(c)$, $n\geq 4$,  be an isometric immersion of a simply connected Riemannian manifold. 
Assume that  
$f$  has a principal curvature $\lambda$ of (constant) multiplicity either $n-1$ or $n$.
Then  $M^n$ admits an isometric immersion $\tilde f\colon M^{n} \to \Q_{\tilde s}^{n+1}(\tilde c)$,
 unless $c>\tilde c$, $s=0$ and $\tilde s=1$,  or $c<\tilde c$, $s=1$ and $\tilde s=0$, and assertions $(i)$-$(iii)$ in Theorem \ref{ngeq4} hold.
 Moreover, $\tilde f$  is unique up to congruence except in case $(ii)$.
\end{theorem}

The next result, proved by  Dajczer and the second author in \cite{dt1} when $s=0=\tilde s$, 
shows how any solution 
$f\colon M^{n} \to \Q_s^{n+1}(c)$, $n\geq 4$, of Problem $*$ arises. 

\begin{theorem}\label{prob2n4}
Let $f\colon M^{n} \to \Q_s^{n+1}(c)$ and $\tilde{f}\colon M^{n} \to \Q^{n+1}_{\tilde s}(\tilde{c})$, $n\geq 4$, be isometric immersions with, say, $c>\tilde c$. 
If  $s=0$, assume that $\tilde s=0$. Then, for  $s=\tilde s$ (respectively,   $s=1$ and $\tilde s=0$),
 there exist, locally on an open dense subset of $M^n$,  isometric embeddings 
 $$H\colon \Q_s^{n+1}(\tilde c)\to \Q^{n+2}_{s}(\tilde{c})\,\,\,\,\mbox{and}\,\,\,\,\,
 i\colon \Q_s^{n+1}(c)\to \Q^{n+2}_{s}(\tilde{c})$$
  (respectively, $H\colon \Q_s^{n+1}(c)\to \Q^{n+2}_{s}({c})$
 and $i\colon \Q_{\tilde s}^{n+1}(\tilde c)\to \Q^{n+2}_{s}({c})$), 
 with $i$  umbilical,    and an isometry 
 $$\Psi\colon \bar M^n:=H(\Q_s^{n+1}(\tilde c))\cap i(\Q_s^{n+1}(c))\to M^n$$
 (respectively, $\Psi\colon \bar M^n:=H(\Q_s^{n+1}(c))\cap i(\Q_{\tilde s}^{n+1}(\tilde c))\to M^n$) 
 such that
$$f\circ \Psi= i^{-1}|_{\bar M^n}\,\,\,\mbox{and}\,\,\,\,\tilde f\circ \Psi= H^{-1}|_{\bar M^n}.$$
 (respectively, 
$f\circ \Psi= H^{-1}|_{\bar M^n}$ and $\tilde f\circ \Psi= i^{-1}|_{\bar M^n}$).
\end{theorem}

Theorem \ref{prob2n4} explains the existence of a principal curvature $\lambda$ of multiplicity 
at least $n-1$ for a solution 
$f\colon M^{n} \to \Q_s^{n+1}(c)$, $n\geq 4$, of Problem $*$ : the (images by $f$ of the) leaves
of the distribution on $M^n$ given by the eigenspaces of $\lambda$ are the intersections with $i(\Q_s^{n+1}(\tilde c))$ of the (images by $H$ of the) relative nullity leaves of $H$, which have dimension at least $n$. \vspace{2ex}

Next we consider Problem $*$  for hypersurfaces of dimension $n=3$. 
The following result provides the solutions in two (``dual") special cases.

\begin{theorem}\label{multdoisb}
Let $f\colon M^{3} \to \Q_s^{4}(c)$ be a hypersurface for which there exists  an isometric immersion   $\tilde{f}\colon M^{3} \to \Q^{4}_{\tilde s}(\tilde{c})$ with  
 $\tilde c \neq  {c}$. 
\begin{itemize}
 \item[(a)]  Assume that  $f$ has a principal curvature of multiplicity two. If either $c> \tilde{c}$,    $s=0$ and $\tilde s=1$, or if  $c < \tilde{c}$, 
  $s=1$ and $\tilde s=0$, then 
$f$ is a rotation hypersurface whose profile curve  is a $\tilde c$-helix in a totally geodesic  surface $\Q_s^{2}(c)$ of $\Q_s^{4}(c)$ and $\tilde f$ is a generalized cone over a surface 
with constant curvature in an umbilical hypersurface $\Q^{3}_{\tilde s}(\bar{c})$ of $\Q^{4}_{\tilde s}(\tilde{c})$, $\bar c\geq \tilde c$. Otherwise, either the same conclusion holds or $f$ and $\tilde f$ are locally given on an open dense subset  as described in Theorem \ref{prob2n4}.
\item[(b)] If one of the principal curvatures of $f$ is zero, then $f$ is a generalized cone over a surface 
with constant curvature in an umbilical hypersurface $\Q^{3}_{s}(\bar{c})$ of $\Q^{4}_{s}({c})$, $\bar c\geq c$,  and $\tilde f$ is a rotation hypersurface whose profile curve 
is a $c$-helix in a totally geodesic  surface $\Q_{\tilde s}^{2}(\tilde c)$ of $\Q_{\tilde s}^{4}(\tilde c)$.
\end{itemize}
 \end{theorem}
 
 By a \emph{generalized cone} over a surface 
$g\colon M^2\to \Q_s^3(\bar c)$ in an umbilical hypersurface $\Q_s^3(\bar c)$ of $\Q^{4}_{s}({c})$, $\bar c\geq c$, we mean the 
hypersurface parametrized by (the restriction to the subset of regular points of) the map $G\colon M^2\times \R\to \Q^{4}_{s}({c})$ given by $$G(x, t)=\exp_{g(x)}(t\xi(g(x))),$$ where $\xi$ is a unit normal vector field to the inclusion 
$i\colon \Q_s^3(\bar c)\to \Q^{4}_{s}({c})$ and $\exp$ is the exponential map of  $\Q^{4}_{s}({c})$. A \emph{$c$-helix} in  $\Q_s^2(\tilde c)\subset \R^3_{s+\epsilon_0}$ with respect to a unit vector $v\in \R^3_{s+\epsilon_0}$ is a unit-speed curve $\gamma\colon I\to \Q_s^2(\tilde c)\subset \R^3_{s+\epsilon_0}$ such that the height function $\gamma_v=\<\gamma, v\>$ satisfies
$\gamma_v''+c\gamma_v=0$. Here $\epsilon_0=0$ or $1$, corresponding to $\tilde c>0$ or $\tilde c<0$, respectively.\vspace{1ex}

 In order to deal with the generic case of Problem $*$ for hypersurfaces of dimension $3$, we need to recall the notion of holonomic hypersurfaces.  We call a hypersurface $f\colon M^{n} \to \Q_s^{n+1}(c)$  \emph{holonomic} if $M^n$ carries  global orthogonal coordinates  $(u_1,\ldots, u_{n})$  such that the coordinate vector fields 
$\d_j=\dfrac{\partial}{\partial u_j}$ are everywhere eigenvectors of the shape operator $A$ of $f$. 
Set $v_j=\|\d_j\|$,
 and define $V_{j} \in C^{\infty}(M)$, $1\leq j\leq n$, by $A\partial_j=v_j^{-1}V_{j}\partial_j.$
Thus, the first and second fundamental forms of $f$ are
\be\label{fundforms}
I=\sum_{i=1}^nv_i^2du_i^2\,\,\,\,\mbox{and}\,\,\,\,\,II=\sum_{i=1}^nV_iv_idu_i^2.
\ee
Set $v=(v_1,\ldots, v_n)$ and  $V=(V_{1},\ldots, V_n)$. We call  $(v,V)$ the pair associated to $f$. 
The next result is well known.

\begin{proposition}\label{fund}
 The triple $(v,h,V)$, where 
 $h_{ij}=\frac{1}{v_i}\frac{\partial v_j}{\partial u_i},$
 satisfies the  system of PDE's
 
  \begin{eqnarray}\label{sistema-hol}
 \left\{\begin{array}{l}
  (i) \dfrac{\partial v_i}{\partial u_j}=h_{ji}v_j,\,\,\,\,\,\,\,\,\,(ii) \dfrac{\partial h_{ik}}{\partial u_j}=h_{ij}h_{jk},\vspace{1ex}\\
  (iii) \dfrac{\partial h_{ij}}{\partial u_i} + \dfrac{\partial h_{ji}}{\partial u_j} + h_{ki}h_{kj} + 
  \epsilon V_{i}V_{j}+cv_iv_j=0,\vspace{1ex}\\
  (iv) \dfrac{\partial V_{i}}{\partial u_j}=h_{ji}V_{j},\,\,\,\,1\leq i \neq j \neq k \neq i\leq n.
  \end{array}\right.
 \end{eqnarray}
Conversely, if $(v,h,V)$ is a solution of (\ref{sistema-hol}) on a simply connected open subset   $U \subset \R^{n}$, with  $v_i \neq 0$
everywhere  for all $1\leq i\leq n$,
then there exists a holonomic hypersurface $f\colon U \to  \Q_s^{n+1}(c)$  whose first and second fundamental forms are given by  (\ref{fundforms}).
\end{proposition}

The following characterization of hypersurfaces $f\colon\,M^3\to\Q_s^{4}(c)$ with three distinct principal curvatures that are solutions of Problem~$*$ is one of the main results of the paper. 

\begin{theorem}\label{main2} Let $f\colon\,M^3\to\Q_s^{4}(c)$ 
be a simply connected  holonomic  hypersurface whose associated pair $(v, V)$ satisfies
\be\label{sphol}\sum_{i=1}^3\delta_iv_i^2=\hat \epsilon, \,\,\,\,\,\,\sum_{i=1}^3\delta_iv_iV_i=0\,\,\,\,\,\mbox{and}\,\,\,\,\,\sum_{i=1}^3\delta_iV_i^2=C:=\tilde \epsilon(c-\tilde c), 
\ee
where $\hat \epsilon, \tilde \epsilon\in \{-1, 1\}$, 
$(\delta_1, \delta_2, \delta_3)=(1, -1, 1)$    either if $\hat \epsilon=1$ or if  $\hat \epsilon=-1$ and 
$C>0$,
and   $(\delta_1, \delta_2, \delta_3)=(-1, -1, -1)$ if  $\hat \epsilon=-1$ and 
$C<0$. Then $M^3$ 
admits an isometric immersion into $\Q_{\tilde s}^{4}(\tilde c)$,  which is unique up to congruence.


Conversely, if 
$f\colon M^{3} \to \Q_s^{4}(c)$ is a hypersurface with three distinct principal curvatures for which  there exists an isometric immersion  $\tilde{f}\colon M^{3} \to \Q^{4}_{\tilde s}(\tilde{c})$ with $\tilde c \neq c$, 
then $f$ is locally a  holonomic hypersurface whose associated pair $(v, V)$ satisfies (\ref{sphol}). 
\end{theorem}



As we shall make precise in the sequel, the class of hypersurfaces that are solutions of  Problem $*$ is
closely related to that  of conformally flat hypersurfaces of $\Q_s^{n+1}(c)$, that is, isometric immersions 
$f\colon M^{n} \to \Q_s^{n+1}(c)$ of conformally flat manifolds. Recall that a Riemannian manifold $M^n$ 
is \emph{conformally flat} if each point of $M^n$ has an open neighborhood that is conformally diffeomorphic to 
an open subset of Euclidean space $\R^n$.
First, for $n\geq 4$ we have the following extension of a result due to  E. Cartan when $s=0$.

\begin{theorem}\label{cfhyp}
Let $f\colon M^{n} \to \Q_s^{n+1}(c)$  be a hypersurface of   dimension   $n\geq 4$. Then $M^n$ is conformally flat
 if and only if  $f$  has a principal curvature of multiplicity 
at least $n-1$.
\end{theorem}

It was already known by E. Cartan that the ``only if" assertion in the preceding result is no longer true for $n=3$ and $s=0$. The study of conformally flat hypersurfaces by Cartan was  taken up by  Hertrich-Jeromin \cite{h-j}, who showed  that a conformally flat hypersurface $f\colon\,M^3\to\Q^{4}(c)$ with three distinct principal curvatures
 admits locally principal coordinates $(u_1, u_2, u_3)$ such that the induced metric $ds^2=\sum_{i=1}^3v_i^2du_i^2$ satisfies, say,
$v_2^2=v_1^2+v_3^2$. The next result states that conformally flat hypersurfaces $f\colon\,M^3\to\Q_s^{4}(c)$ with three distinct 
principal curvatures are characterized by the existence of such principal coordinates
under some additional conditions. 

\begin{theorem}\label{main3} Let $f\colon\,M^3\to\Q_s^{4}(c)$ 
be a  holonomic hypersurface whose associated pair $(v, V)$ satisfies 
\be\label{eq:cflatcond}\sum_{i=1}^3\delta_iv_i^2=0, \,\,\,\,\,\,\sum_{i=1}^3\delta_iv_iV_i=0\,\,\,\,\,\mbox{and}\,\,\,\,\,\sum_{i=1}^3\delta_iV_i^2=1, 
\ee
where $(\delta_1, \delta_2, \delta_3)= (1,-1, 1)$. Then $M^3$ is conformally flat. 

Conversely, any conformally flat hypersurface $f\colon\,M^3\to\Q_s^{4}(c)$ with three distinct 
principal curvatures is locally a 
holonomic hypersurface whose associated pair $(v, V)$ satisfies (\ref{eq:cflatcond}).
\end{theorem}

It is amazing that the class of holonomic Euclidean hypersurfaces of any dimension $n$  whose associated pair $(v, V)$ satisfies the conditions 
$$\sum_{i=1}^n\delta_iv_i^2=K_1\,\,\,\,\,\mbox{and}\,\,\,\,\,\sum_{i=1}^n\delta_iV_i^2=K_2,$$
where  $K_1, K_2\in \R$ and $\delta_i\in \{-1, 1\}$ for $1\leq i\leq n$, was considered by Bianchi \cite{bi} almost one century ago,
his interest on such hypersurfaces relying on the fact that they satisfy many of the properties of  constant curvature surfaces an their parallel surfaces in  $\R^3$. 
In particular, a Ribaucour transformation  for that class was sketched in Bianchi's paper.

It follows from Theorems \ref{main2} and \ref{main3} that, in order to produce  hypersurfaces of $\Q_s^4(c)$ that are either conformally flat or 
admit an isometric immersion into  $\Q_s^4(\tilde c)$ with $\tilde c\neq c$, one must start with 
solutions $(v,h,V)$ on an open simply connected subset $U\subset \R^3$ of the same system of PDE's, namely, the one obtained by adding to system (\ref{sistema-hol}) (for $n=3$) the equations
\be\label{eq:int1}
\delta_i\dfrac{\partial v_i}{\partial u_i}+ \delta_j h_{ij}v_j+ \delta_k h_{ik}v_k=0
\ee
and
\be\label{eq:int2}
\delta_i\dfrac{\partial V_i}{\partial u_i}+ \delta_j h_{ij}V_j+ \delta_k h_{ik}V_k=0, \,\,\,\,1\leq i\neq j\neq k\neq i\leq 3,
\ee
with $(\delta_1, \delta_2, \delta_3)=(1, -1, 1)$. Such system has the first integrals
$$\sum_{i=1}^3\delta_iv_i^2=K_1,\,\,\,\,\,\,\,\sum_{i=1}^3\delta_iv_iV_i=K_2\,\,\,\,\mbox{and}\,\,\,\,\,\sum_{i=1}^3\delta_iV_i^2=K_3.$$
If initial conditions at some point are chosen so that $K_1=1$ (respectively, $K_1=0$), $K_2=0$ and $K_3=\epsilon(c-\tilde c)$ 
(respectively, $K_3=1$), then the corresponding solutions give rise to hypersurfaces of $\Q_s^4(c)$ with three distinct principal curvatures that can be isometrically immersed into $\Q_s^{4}(\tilde c)$
(respectively, are conformally flat).

  Our characterizations in Theorems \ref{main2} and \ref{main3} of  hypersurfaces  of  $\Q_s^4(c)$ with three distinct principal curvatures that admit an isometric immersion 
  into  $\Q_{\tilde s}^4(\tilde c)$, with $c\neq \tilde c$, or are conformally flat, respectively, allow us to derive a Riabaucour transformation for both classes of hypersurfaces.
  In particular, it yields the following 
  process to generate a family of new elements of such classes 
   from a given one. 
  We denote by $i\colon \Q_s^{4}(c)\to \R_{s+\epsilon_0}^{5}$  an umbilical inclusion,  where $\epsilon_0=0$ or $1$, corresponding to $c>0$ or $c<0$, respectively.

  \begin{theorem}\label{ribA}
If  $f \colon M^3 \to \Q_s^4(c)$ is a  holonomic  hypersurface whose associated pair $(v,  V)$ satisfies  (\ref{sphol}) (respectively,  (\ref{eq:cflatcond})), then the linear system of PDE's 
\begin{eqnarray}\label{sist-RibA}
\left\{\begin{array}{l}
(i) \dfrac{\partial \varphi}{\partial u_i} = v_i\gamma_i,\,\,\,\,\,\,\,
(ii) \dfrac{\partial \gamma_j}{\partial u_i} = h_{ji}\gamma_i,\ \ \,\,\ \ i\neq j,\vspace{1ex}\\
(iii) \dfrac{\partial \gamma_i}{\partial u_i} = (v_i - v_i')\psi - \sum_{j \neq i}h_{ji}\gamma_j + \beta V_i - c \varphi v_i, \vspace{1ex}\\
(iv) \epsilon\dfrac{\partial \beta}{\partial u_i} = -V_{i}\gamma_i,\vspace{1ex}\\
(v) \dfrac{\partial \log \psi}{\partial u_i} = -\dfrac{\gamma_i v_i'}{\varphi},\,\,\,\,\,\,\,\,
(vi) \dfrac{\partial v_i'}{\partial u_j} = h'_{ji}v_j',\ \ \,\,\ \ i\neq j,\vspace{1ex}\\
(vii) \delta_i\dfrac{\partial v_i'}{\partial u_i}+ \delta_j h'_{ij}v_j'+ \delta_k h'_{ik}v_k'=0,
\end{array}\right.
\end{eqnarray}
%
where 
\be\label{eq:hij'}
h_{ij}=\frac{1}{v_i}\frac{\partial v_j}{\partial u_i}\,\,\,\,\,\mbox{and}\,\,\,\,\,h'_{ij}=h_{ij} + (v'_j - v_j)\dfrac{\gamma_i}{\varphi},
\ee
 is completely integrable and has the first integrals 
\begin{equation}\label{int-primeiraA}
\sum_i {\gamma_i}^2+\epsilon\beta^2 + c\varphi^2 - 2\varphi \psi = K_1\in \R
\end{equation}
and 
\begin{equation}\label{int-primeiraB}
\delta_1 v_1'^2 + \delta_2 v_2'^2 + \delta_3 v_3'^2 =  K_2\in \R. 
\end{equation}
Let $(\gamma_1,\gamma_2,\gamma_3,v_1',v_2',v_3',\varphi,\psi,\beta)$ be a solution of (\ref{sist-RibA}) with initial conditions at some point chosen so that $K_1=0$ and $K_2=\hat \epsilon$ (respectively, $K_2=0$), and so that   the  function 
\be\label{eq:omega}
\Omega=\varphi\sum_{j=1}^3\delta_j v_j'V_j  - \epsilon\beta \left(K_2 -\sum_{j=1}^3\delta_j  v_jv_j'\right),
\ee
with $K_2=\hat \epsilon$ (respectively, $K_2=0$), vanishes at that point.
Then, the map $F'\colon M^3\to \R^{5}_{s+\epsilon_0}$, given in terms of  $F=i\circ f$  by
\begin{equation}\label{dt2:10A}
{F'} = F - \dfrac{1}{\psi}\left(\sum_i \gamma_i F_*e_i + \beta i_*\xi + c\varphi F \right),
\end{equation}
where $\xi$ is a  unit normal vector field $\xi$ to $f$ and $e_i=v_i^{-1}\d_i$, $1\leq i\leq 3$,
satisfies $F'=i\circ  f'$, where $f'\colon M^3\to \Q_s^4(c)$ is a  holonomic  hypersurface whose associated pair $(v',  V')$, with
$$V_i' = V_i + (v_i - v_i')\dfrac{\epsilon\beta}{\varphi},$$
 also satisfies  (\ref{sphol}) (respectively,  (\ref{eq:cflatcond})).
\end{theorem}


Explicit examples of   hypersurfaces  of  $\Q_s^4(c)$ with three distinct principal curvatures that admit an isometric immersion 
  into  $\Q_{\tilde s}^4(\tilde c)$ with $c\neq \tilde c$, as well as of  conformally flat hypersurfaces  of  $\R_s^4$ with three distinct principal curvatures, are constructed in Section $6$ by means of Theorem \ref{ribA}.

As a special consequence of Theorem \ref{ribA}, it follows that   hypersurfaces $f \colon M^3 \to \Q_s^4(c)$ that can be  isometrically immersed
into $\R_{\tilde s}^4$ 
arise in families of parallel hypersurfaces.

\begin{corollary}\label{parallel}
Let  $f\colon M^3 \to \Q_s^4(c)$ be a holonomic  hypersurface whose associated pair $(v,  V)$ satisfies  (\ref{sphol}) with $\tilde c=0$.  Then any parallel hypersurface
$f_t\colon M^3 \to \Q_s^4(c)$ to $f$ has also the same property. 
\end{corollary}

It was already shown in \cite{dt1} for $s=0=\tilde s$  that, unlike the case of dimension $n\geq 4$, among hypersurfaces 
$f\colon M^n\to \Q_s^{n+1}(c)$ of  dimension $n=3$ with three distinct 
principal curvatures,  the classes of  solutions  of Problem $*$ and conformally flat  hypersurfaces 
are distinct. Moreover,  it was observed that their intersection contains the generalized cones over surfaces 
with constant curvature in an umbilical hypersurface $\Q_s^{3}(\bar c)$ of $\Q^{4}_{s}({c})$, $\bar c\geq c$.
 Our last result states that such intersection contains no other elements.

\begin{theorem}\label{intersection} Let $f\colon M^3\to \Q_s^{4}(c)$ be a conformally flat hypersurface with three distinct principal curvatures.
If $M^3$ admits an isometric immersion into $\Q_{\tilde s}^{4}(\tilde c)$,   $\tilde c\neq  c$, then 
$f$ is a generalized cone over a surface 
with constant curvature in an umbilical hypersurface $\Q_s^{3}(\bar c)$ of $\Q^{4}_{s}({c})$, $\bar c\geq c$.
\end{theorem}

\section{Proofs of Theorems \ref{ngeq4},  \ref{conv} and \ref{prob2n4}}
\noindent \emph{Proof of Theorem \ref{ngeq4}:} Let $i\colon \Q_s^{n+1}(c)\to \Q_{s+\epsilon_0}^{n+2}(\tilde c)$ be an umbilical inclusion, 
where $\epsilon_0=0$ or $1$, corresponding to $c>\tilde c$ or $c<\tilde c$, respectively, and set $\hat f=i\circ f$. 
Then, the second fundamental forms $\alpha$ and $\hat \alpha$ of $f$ and $\hat f$, respectively, are related by
\be\label{sffs} \hat \alpha=i_*\alpha+\sqrt{|c-\tilde c|}\<\,\,,\,\,\>\xi,
\ee
where $\xi$ is one of the unit  vector fields that are normal to $i$.

For a fixed point $x\in M^n$, define  $W^{3}(x):= N_{\hat{f}}M(x) \oplus N_{\tilde{f}}M(x)$, and endow $W^3(x)$ with the inner product 
$$\<\!\<(\xi+\tilde{\xi}, \eta + \tilde{\eta})\>\!\>_{W^3(x)}:= \<\xi,\eta\>_{N_{\hat{f}}M(x)} -\<\tilde{\xi},\tilde{\eta}\>_{N_{\tilde{f}}M(x)},$$ 
which has index $(s+\epsilon_0)+(1-\tilde s)$.                                                                                                                                                                                                        

Now define a bilinear form  $\beta_x\colon T_xM \times T_xM \to W^3(x)$ by $$\beta_x=\hat{\alpha}(x) \oplus \tilde{\alpha}(x),$$ where $\hat\alpha(x)$ and $\tilde{\alpha}(x)$ 
are the second fundamental forms of $\hat f$ and $\tilde f$, respectively,  at $x$.  Notice that ${\cal N}(\beta_x)\subset {\cal N}(\hat \alpha(x))=\{0\}$ by (\ref{sffs}). 
On the other hand, it follows from the Gauss equations of  $\hat{f}$  and  $\tilde{f}$ that $\beta_x$ is flat with respect to $\<\!\<\,\,,\,\,\>\!\>$, that is,
$$\<\!\<\beta_x(X,Y),\beta_x(Z,W)\>\!\>=\<\!\<\beta_x(X,W),\beta_x(Z,Y)\>\!\>$$
for all $X,Y,Z, W\in T_xM$. 
Thus, if $\<\!\<\,\,,\,\,\>\!\>$ is positive definite, which is the case when $s=0$, $\tilde s=1$ and $\epsilon_0=0$, that is, $c>\tilde c$, we obtain a contradiction with Corollary $1$ of \cite{mo},
according to which one has the inequality
\be\label{moore}\dim {\cal N}(\beta_x)\geq n-\dim W(x)=n-3>0.\ee
The same contradiction is reached by applying the preceding inequality to $-\<\!\<\,\,,\,\,\>\!\>$ when $s=1$, $\tilde s=0$ and $c<\tilde c$, in which case $\<\!\<\,\,,\,\,\>\!\>$ is negative definite.
Therefore, such cases can not occur, which proves the first assertion. 

In all other cases, the index of  $\<\!\<\,\,,\,\,\>\!\>$ is either $1$ or $2$. Thus, by applying Corollary $2$ in \cite{mo} to $\<\!\<\,\,,\,\,\>\!\>$ in the first case and to 
$-\<\!\<\,\,,\,\,\>\!\>$ in the latter, we obtain that ${\cal S}(\beta_x)$ must be degenerate, for otherwise the inequality (\ref{moore}) would still hold, and then we would reach a contradiction 
as before. 

Since ${\cal S}(\beta_x)$ is degenerate, there exist  $\zeta \in N_{\hat{f}}M(x)$ and  $\tilde{N} \in N_{\tilde{f}}M(x)$ such that $(0,0)\neq (\zeta, \tilde{N})\in 
{\cal S}(\beta_x)\cap {\cal S}(\beta_x)^\perp$. In particular, from 
$0=\<\!\<\zeta + \tilde{N},\zeta+ \tilde{N}\>\!\>$ it follows that
$\<\tilde N, \tilde N\>=\<\zeta, \zeta\>$. Thus, either $\tilde N=0$ and $\zeta\in {\cal S}(\hat \alpha(x))\cap {\cal S}(\hat \alpha(x))^\perp$, or we can assume that   $\<\tilde N, \tilde N\>=\tilde \epsilon=\<\zeta, \zeta\>$.

The former case occurs precisely when $f$ is umbilical at $x$ with a principal curvature $\lambda$ with respect to one of the unit normal vectors $N$ to $f$, satisfying 
$$\epsilon\lambda^2+c-\tilde c=0,$$
in which case  $N_{\hat f}M(x)$ is  a Lorentzian two-plane and  $\zeta=\lambda i_*N+\sqrt{|c-\tilde c|}\xi$ is a light-like vector that spans ${\cal S}(\hat \alpha(x))$. 
In this case,  all sectional curvatures of $M^n$ at $x$ are equal to $\tilde c$ by the Gauss equation of $f$, and hence $\tilde f$    has $0$ as a principal curvature at $x$ with  multiplicity at least $n-1$ by the Gauss equation of $\tilde f$. 

Now assume that $\<\tilde N, \tilde N\>=\tilde \epsilon=\<\zeta, \zeta\>$. Then, from
$$
0=\<\!\<\beta, \zeta + \tilde{N} \>\!\> =  \<\hat{\alpha},\zeta\> - \<\tilde{\alpha},\tilde{N}\>,
$$
we obtain that $A^{\hat f}_{\zeta} = A^{\tilde f}_{\tilde{N}}$. Let $\zeta^\perp\in N_{\hat f}M(x)$ be such that $\{\zeta, \zeta^\perp\}$ is an orthonormal basis of  $N_{\hat f}M(x)$. 
The Gauss equations for $\hat f$ and $\tilde f$ imply that 
$$\<A^{\hat f}_{\zeta^\perp}X,Y\>\<A^{\hat f}_{\zeta^\perp}Z,W\>=\<A^{\hat f}_{\zeta^\perp}X,W\>\<A^{\hat f}_{\zeta^\perp}Z,Y\>$$
for all $X, Y, Z, W\in T_xM$, which is equivalent to $\dim {\cal N}(A^{\hat f}_{\zeta^\perp})\geq n-1$. 
Since $A^{\hat f}_{\xi}=\delta\sqrt{|c-\tilde c|}I$ by (\ref{sffs}), with $\delta=(c-\tilde c)/|c-\tilde c|$, it follows that the restriction to
${\cal N}(A^{\hat{f}}_{\zeta^{\perp}})$ of all shape operators 
$A^{\hat f}_{\eta}$, $\eta\in N_{\hat f}M(x)$, is a multiple of the identity tensor. In particular, this is the case for 
$A^{\hat f}_{i_*N}=A^f_N$, where $N$ is one of the unit normal vector fields to $f$, hence $f$ has a principal curvature $\lambda$ at $x$ with 
multiplicity at least $n-1$. 

Moreover, if $\lambda=0$ then $\zeta^\perp$ must coincide with $i_*N$, and hence $\zeta$ with $\xi$, up to signs. 
Therefore $A^{\tilde f}_{\tilde N}=A^{\hat f}_{\xi}$, up to sign,  hence $\tilde f$ is umbilical at $x$. 
If $f$ is umbilical at $x$ and $c+\epsilon\lambda^2\neq \tilde c$, then $A_{\zeta^\perp}=0$ and 
$A^{\tilde f}_{\tilde N}=A^{\hat f}_{\zeta}$ is a (nonzero) constant multiple of the identity tensor. 
Finally, if $\lambda\neq 0$ has multiplicity $n-1$, then we must have $\zeta^\perp\neq i_*N$ and $\dim {\cal N}(A^{\hat f}_{\zeta^\perp})= n-1$, hence ${\cal N}(A^{\hat f}_{\zeta^\perp})$ is an eigenspace of $A^{\hat f}_{\zeta}=A^{\tilde f}_{\tilde N}$.
\vspace{2ex}\qed

\noindent \emph{Proof of Theorem \ref{conv}:}  Suppose first that $f$ is umbilical, with a (constant) principal curvature $\lambda$.
 If $c+\epsilon\lambda^2= \tilde c$, then $M^n$ has constant curvature $\tilde c$, hence it  admits isometric immersions into $\Q_{\tilde s}^{n+1}(\tilde c)$ having $0$ as a principal curvature with multiplicity at least $n-1$. Otherwise, by the assumption there exists $\tilde \lambda\neq 0$ such that 
$c-\tilde c+\epsilon \lambda^2=\tilde \epsilon \tilde \lambda^2$.
Hence
$c+\epsilon\lambda^2=\tilde c+\tilde \epsilon\tilde \lambda^2,$
thus $\tilde A=\tilde \lambda I$ satisfies the Gauss and Codazzi equation for an (umbilical) isometric immersion into 
$\Q_{\tilde s}^{n+1}(\tilde c)$.

Assume  now that  
$f$  has  principal curvatures $\lambda$ and $\mu$  of multiplicities  $n-1$ and $1$, respectively, with corresponding eigenbundles $E_\lambda$ and $E_\mu$. If $\lambda=0$, then $M^n$ has constant curvature $c$, hence it admits an umbilical isometric immersion into 
$\Q_{\tilde s}^{n+1}(\tilde c)$. From now on, assume that $\lambda\neq 0$. 
Then, one can check that  the Codazzi equations for $f$ are equivalent to the fact that $E_\lambda$ and $E_\mu$ are umbilical distributions with mean curvature normals $\eta$ and $\zeta$,
respectively, satisfying
$$\eta=\frac{(\nabla \lambda)_{E_\mu}}{\lambda -\mu}\,\,\,\,\mbox{and}\,\,\,\,\zeta=\frac{(\nabla \mu)_{E_\lambda}}{\mu-\lambda}.$$

By the assumption, there exist   $\tilde \lambda, \tilde \mu\in C^{\infty}(M)$  such that 
$$c-\tilde c+\epsilon \lambda^2=\tilde \epsilon \tilde \lambda^2 \,\,\,\,\,\mbox{and}\,\,\,\,\,c-\tilde c+\epsilon \lambda\mu=\tilde \epsilon \tilde \lambda\tilde \mu.$$
Moreover, the first of the preceding equations implies that $\tilde \lambda\neq 0$ everywhere, and hence $\tilde \lambda$ and $\tilde \mu$ are unique if $\tilde \lambda$ is chosen to be positive.
From both equations we obtain that
$$\epsilon \lambda^2-\tilde \epsilon \tilde \lambda^2=\epsilon \lambda\mu-\tilde \epsilon \tilde \lambda\tilde \mu,\,\,\,\,\,\,\,\,\,
\epsilon \lambda\nabla\lambda=\tilde\epsilon \tilde\lambda\nabla\tilde\lambda$$
and
$$\epsilon ((\nabla\lambda)\mu+\lambda\nabla\mu)=\tilde\epsilon ((\nabla\tilde\lambda)\tilde\mu+ \tilde\lambda\nabla\tilde\mu).$$
It follows that
\be\label{cod1}
\frac{(\nabla \tilde\lambda)_{E_\mu}}{\tilde\lambda -\tilde\mu}=\frac{(\nabla \lambda)_{E_\mu}}{\lambda -\mu}
\ee
and similarly,
\be\label{cod2}
\frac{(\nabla \tilde\mu)_{E_\lambda}}{\tilde\mu-\tilde\lambda}=\frac{(\nabla \mu)_{E_\lambda}}{\mu-\lambda}.
\ee
Let $\tilde A$ be the endomorphism of $TM$ with eigenvalues $\tilde \lambda$, $\tilde \mu$ and corresponding eigenbundles $E_\lambda$ and $E_\mu$, respectively.
Since 
$$c+\epsilon\lambda^2=\tilde c+\tilde \epsilon\tilde \lambda^2\,\,\,\,\,\mbox{and}\,\,\,\,c+\epsilon\lambda\mu=\tilde c+\tilde \epsilon \tilde\lambda\tilde \mu,$$
the Gauss equations for an isometric immersion $\tilde{f}\colon M^{n} \to \Q^{n+1}_{\tilde s}(\tilde{c})$ are satisfied by $\tilde A$. It follows from 
(\ref{cod1}) and (\ref{cod2}) that $\tilde A$ also satisfies the Codazzi equations.\vspace{2ex}\qed\\
\noindent \emph{Proof of Theorem \ref{prob2n4}:} Since we are assuming that $c>\tilde c$, there exist  umbilical inclusions 
$i\colon \Q_s^{n+1}(c)\to \Q^{n+2}_{s}(\tilde{c})$ and  $i\colon \Q_{\tilde s}^{n+1}(\tilde c)\to \Q^{n+2}_{s}({c})$ for $(s,\tilde s)=(1,0)$. If $s=\tilde s$ (respectively, $(s,\tilde s)=(1,0)$), set  $\hat f=i\circ f$ (respectively, $\hat f=i\circ \tilde f$). 
Then, one can use the existence of normal vector fields $\zeta\in \Gamma(N_{\hat f}M)$ and $\tilde N\in \Gamma(N_{\tilde f}M)$ satisfying $\<\zeta, \zeta\>=\tilde \epsilon=\<\tilde N, \tilde N\>$ and $A_{\zeta}^{\hat f}=A_{\tilde N}^{\tilde f}$ and argue as in the proof of Theorem $3$ in \cite{dt1}. One obtains that there exists an open dense subset $U\subset M^n$, each point of which has an open neighborhood $V\subset M^n$ such that $\hat f|_V$ (respectively, $f|_V$) is a composition $\hat f|_V=H\circ \tilde f|_V$ (respectively, $f|_V=H\circ \hat f|_V$) with an isometric embedding $H\colon W\subset \Q_s^{n+1}(\tilde c)\to \Q^{n+2}_{s}(\tilde{c})$ (respectively,  $H\colon W\subset \Q_s^{n+1}(c)\to \Q^{n+2}_{s}({c})$), with $\tilde f(V)\subset W$ (respectively, $\hat f(V)\subset W$). Set $\bar M^n=H(W)\cap i(\Q_s^{n+1}(c))$ (respectively, 
$\bar M^n=H(W)\cap i(\Q_{\tilde s}^{n+1}(\tilde c))$). Then $i\circ f|_V=H\circ \tilde f|_V\colon V\to \bar M^n$ (respectively, 
$H\circ f|_V=i\circ \tilde f|_V\colon V\to \bar M^n$) is an isometry. Let $\Psi\colon \bar M^n\to V$ be the inverse of this isometry. Then $f\circ \Psi= i^{-1}|_{\bar M^n}$ and $\tilde f\circ \Psi= H^{-1}|_{\bar M^n}$ (respectively, 
$f\circ \Psi= H^{-1}|_{\bar M^n}$ and $\tilde f\circ \Psi= i^{-1}|_{\bar M^n}$), where $i^{-1}$ and $H^{-1}$ denote the inverses of the maps $i$ and $H$, respectively, regarded as maps onto their images. \qed

\section{Proof of Theorem \ref{multdoisb}}

Before going into the proof of Theorem \ref{multdoisb}, we establish a basic fact that will also be used in the proof of Theorem \ref{main2} in the next section.

\begin{lemma}\label{base-comum}
Let $f\colon M^{3} \to \Q_s^{4}(c)$ and  $\tilde{f}\colon M^{3} \to \Q^{4}_{\tilde s}(\tilde{c})$ be hypersurfaces with  $c \neq \tilde{c}$. Then, at each point $x\in M^3$ there exists an  orthonormal basis $\{e_1,e_2,e_{3}\}$ of $T_xM^{3}$ that simultaneously  diagonalizes the second fundamental forms of   $f$ and $\tilde{f}$. 
 \end{lemma}
 \proof Define  $i\colon \Q_s^{4}(c)\to \Q_{s+\epsilon_0}^{5}(\tilde c)$ and  $\hat f$, as well as
  $W^{3}(x)$,  $\<\!\<\,\,, \,\,\>\!\>_{W^3(x)}$ and   $\beta_x$ for each $x\in M^n$,   as in the proof of Theorem \ref{ngeq4}. 
  If ${\cal S}(\beta_x)$ is degenerate for all $x\in M^3$, we conclude as in the case $n\geq 4$ that the assertions in Theorem~\ref{ngeq4} hold, hence 
  the statement is clearly true in this case. 
  
  Suppose now that ${\cal S}(\beta_ x)$ is nondegenerate at $x\in M^3$. Then the inequality 
  $$\dim {\cal S}(\beta_x)\geq \dim T_xM-\dim {\cal N}(\beta_x)$$
  holds by Corollary $2$ in \cite{mo}. Since ${\cal N}(\beta_x)=\{0\}$, the right-hand-side is equal to $\dim T_xM=3= \dim W^3(x)$, hence we must have equality  in the above inequality. 
  By Theorem $2-b$ in \cite{mo}, there exists an orthonormal basis  $\{\xi_1,\xi_2,\xi_{3}\}$ of $W^3(x)$ and a  basis $\{\theta^1,\theta^2,\theta^{3}\}$ of $T_x^*M$ such that 
\begin{eqnarray*}
 \beta=\sum_{j=1}^{3}\theta^j \otimes \theta^j \xi_j.
\end{eqnarray*}
In particular, if $i \neq j$ then $\beta(e_i,e_j)=0$  for the dual basis  $\{e_1,e_2,e_{3}\}$ of  $\{\theta^1,\theta^2,\theta^{3}\}$. 
It follows that  $\{e_1,e_2,e_{3}\}$ diagonalyzes both $\hat \alpha$ and $\tilde \alpha$, and therefore both $\alpha$ and $\tilde \alpha$, in view of 
(\ref{sffs}). It also follows from (\ref{sffs}) that 
$$0=\<\hat \alpha(e_i, e_j), \xi\>=\sqrt{|c-\tilde c|}\<e_i, e_j\>,\,\,\,i\neq j,$$
hence the basis $\{e_1,e_2,e_{3}\}$  is orthogonal. \vspace{1ex}\qed

 \begin{lemma}\label{multdois}
 Under the assumptions of Lemma \ref{base-comum},  let  $\lambda_1, \lambda_2, \lambda_3$ and $\mu_1, \mu_2, \mu_3$ be the principal curvatures of  $f$ and $\tilde f$ correspondent to  $e_1, e_2$ and $e_3$, respectively. 
 \begin{itemize}
 \item[(a)]  Assume that $f$ has a principal curvature of multiplicity two,  say, that $\lambda_1=\lambda_2:=\lambda$. If either $c> \tilde{c}$,    $s=0$ and $\tilde s=1$, or   $c < \tilde{c}$, 
  $s=1$ and $\tilde s=0$, then  
 $$c-\tilde c+\epsilon \lambda\lambda_3=0,\,\,\,\,\mu_3=0\,\,\,\mbox{and}\,\,\,
c-\tilde c+\epsilon \lambda^2=\tilde \epsilon \mu_1\mu_2.$$
 Otherwise, either the same conclusion holds or 
 $$\mu_1=\mu_2:=\mu,\,\,\,  
 c-\tilde c+\epsilon \lambda^2=\tilde \epsilon \mu^2\,\,\,\mbox{and}\,\,\, c-\tilde c+\epsilon \lambda\lambda_3=\tilde \epsilon \mu\mu_3.$$
 \item[(b)] Assume, say, that $\lambda_3=0$. Then $\mu_1=\mu_2:=\mu$,
 \be\label{e1}
 c-\tilde c+\epsilon \lambda_1\lambda_2=\tilde \epsilon \mu^2
 \ee
 and
 \be\label{e2}
 c-\tilde c=\tilde \epsilon \mu\mu_3.
 \ee
\end{itemize}
 \end{lemma}
 \proof By the Gauss equations for $f$ and $\tilde f$, we have
 \be\label{gftildef} 
 c+\epsilon \lambda_i\lambda_j=\tilde c+\tilde \epsilon \mu_i\mu_j, \,\,\,1\leq i\neq j\leq 3.
 \ee
 $(a)$ If $\lambda_1=\lambda_2:=\lambda$, then the preceding equations are
 \be\label{c1}
 c+\epsilon \lambda^2=\tilde c+\tilde \epsilon \mu_1\mu_2,
 \ee
 \be\label{c2}
 c+\epsilon \lambda\lambda_3=\tilde c+\tilde \epsilon \mu_1\mu_3
 \ee
 and 
 \be\label{c3}
 c+\epsilon \lambda\lambda_3=\tilde c+\tilde \epsilon \mu_2\mu_3.
 \ee
 The two last equations yield
 $$\mu_3(\mu_1-\mu_2)=0,$$
 hence either $\mu_3=0$ or $\mu_1=\mu_2$. In view of (\ref{c1}), the second possibility can not occur if 
 either  $c> \tilde{c}$,  $s=0$ and $\tilde s=1$, 
 or  $c < \tilde{c}$,  $s=1$ and $\tilde s=0$. Thus, in these cases we must have $\mu_3=0$, and then 
 $c-\tilde c+\epsilon \lambda^2=\tilde \epsilon \mu_1\mu_2$ and $c-\tilde c+\epsilon \lambda\lambda_3=0$ by 
 (\ref{c2}) and (\ref{c3}).
 
 Otherwise, either the same conclusion holds or 
 $\mu_1=\mu_2:=\mu$, and then  
 $c-\tilde c+\epsilon \lambda^2=\tilde \epsilon \mu^2$ and $c-\tilde c+\epsilon \lambda\lambda_3=\tilde \epsilon \mu\mu_3$
 by (\ref{c2}) and (\ref{c3}). 
 \vspace{1ex}\\
 $(b)$ If $\lambda_3=0$, then equations (\ref{gftildef}) become
 \be\label{d1}
 c-\tilde c+\epsilon \lambda_1\lambda_2=\tilde \epsilon \mu_1\mu_2,
 \ee
 \be\label{d2}
 c-\tilde c=\tilde \epsilon \mu_1\mu_3
 \ee
 and 
 \be\label{d3}
c-\tilde c=\tilde \epsilon \mu_2\mu_3
 \ee
 Since $\mu_3\neq 0$ by (\ref{d2}) or (\ref{d3}), these equations imply that  $\mu_1=\mu_2:=\mu$, and we obtain
  (\ref{e2}). Equation (\ref{e1}) then follows from (\ref{d1}).  \vspace{2ex}\qed

\noindent \emph{Proof of Theorem \ref{multdoisb}:}
Assume that  $f$ has a principal curvature of multiplicity two, say, $\lambda_1=\lambda_2:=\lambda$.
Suppose first that  either $c> \tilde{c}$,    $s=0$ and $\tilde s=1$, or   $c < \tilde{c}$, 
  $s=1$ and $\tilde s=0$. Then,  it follows from Lemma \ref{multdois} that 
  \be\label{eq:pcurv}
  c-\tilde c+\epsilon \lambda\lambda_3=0,\,\,\,\,\mu_3=0\,\,\,\mbox{and}\,\,\,
c-\tilde c+\epsilon \lambda^2=\tilde \epsilon \mu_1\mu_2.\ee
In particular, we must have $\lambda\neq 0$ by the first of the preceding equations, whereas the last one  implies that $\mu_1\mu_2\neq 0$. Then, it is well known that  $E_\lambda$ is a spherical distribution,  that is, it is umbilical and its mean curvature normal $\eta=\nu e_3$ satisfies $e_1(\nu)=0=e_2(\nu)$. In particular, a leaf $\sigma$ of $E_\lambda$ has constant sectional curvature $\nu^2+\epsilon\lambda^2+c=\nu^2+\tilde \epsilon\mu_1\mu_2+\tilde c$.  Denoting by $\nabla$ and $\tilde \nabla$  the connections on $M^3$ and $\tilde f^*T\Q_{\tilde s}^{4}(\tilde c)$, respectively, we have
$$\tilde \nabla_{e_i}\tilde{f}_*e_3=\tilde{f}_*\nabla_{e_i}e_3=-\nu \tilde{f}_*e_i,\,\,\,\,1\leq i\leq 2,$$
hence $\tilde{f}(\sigma)$ is contained in an umbilical hypersurface $\Q_{\tilde s}^{3}(\bar c)$ of $\Q_{\tilde s}^{4}(\tilde c)$ with 
constant curvature $\bar c=\tilde c+ \nu^2$ and $\tilde{f}_*e_3$ as a unit normal vector field.

Moreover,  $E_\lambda^\perp=E_{\mu_3}$ is the relative nullity distribution of $\tilde f$. Thus, it is totally geodesic, and in fact its integral curves  are mapped by $\tilde f$ into geodesics of $\Q_{\tilde s}^{4}(\tilde c)$.  It follows that $\tilde f(M^3)$ is contained in a generalized cone over $\tilde f(\sigma)$.   

On the other hand, it is not hard to extend the proof of Theorem $4.2$ in \cite{dcd2} to the case of Lorentzian ambient space forms, and
 conclude that $f$ is a rotation hypersurface in $\Q_s^4(c)$. This means that  there exist subspaces $P^2\subset P^3=P^3_{s+\epsilon_0}$ in $\R^{5}_{s+\epsilon_0}\supset \Q_s^4(c)$ with $P^3\cap \Q_s^4(c)\neq \emptyset$, where $\epsilon_0=0$ or $\epsilon_0=1$, corresponding to $c>0$ or $c<0$, respectively, and a regular curve $\gamma$ in $\Q_s^2(c)=P^3\cap \Q_s^4(c)$ that does not meet $P^2$, such that $f(M^2)$ is the union of the orbits of points of $\gamma$ under the action of the subgroup of orthogonal transformations of $\R^{5}_{s+\epsilon_0}$ that fix pointwise $P^2$. 
 If $P^2$ is nondegenerate, then $f$ can be parameterized by 
 $$
 f(s,u)=(\gamma_1(s)\phi_1(u), \gamma_1(s)\phi_2(u), \gamma_1(s)\phi_3(u),\gamma_{4}(s),\gamma_{5}(s)),
 $$
 with respect to an orthonormal basis $\{e_1, \ldots, e_5\}$ of $\R^5_{s+\epsilon_0}$ satisfying the conditions in either $(i)$ or $(ii)$ below, according to whether the induced metric on $P^2$ has index  
 $s+\epsilon_0$ or $s+\epsilon_0-1$, respectively:
\begin{itemize}
\item[(i)]  $\<e_i, e_i\>=1$ for $1\leq i\leq 3$, $\<e_{3+j}, e_{3+j}\>=\epsilon_j$ for $1\leq j\leq 2$, and $(\epsilon_1, \epsilon_2)$ equal to either $(1,1)$, $(1, -1)$ or $(-1, -1)$, corresponding to 
$s+\epsilon_0=0$, $1$ or $2$, respectively.   
\item[(ii)]  $\<e_1, e_1\>=-1$,  $\<e_i, e_i\>=1$  for $2\leq i\leq 4$  and $\<e_5, e_5\>=\bar{\epsilon}$, where $\bar{\epsilon}=1$ or $\bar{\epsilon}=-1$, corresponding to  $s+\epsilon_0=1$ or $2$, respectively.
\end{itemize}
In both cases,   we have $P^2=\spa\{e_4, e_5\}$,  $P^3=\spa\{e_1, e_4, e_5\}$, $u=(u_1, u_2)$,   $\gamma(s)=(\gamma_1(s), \gamma_4(s), \gamma_5(s))$  a unit--speed curve in $\Q_s^2(c)\subset P^3$ and
$\phi(u)=(\phi_1(u), \phi_2(u), \phi_3(u))$  an orthogonal parameterization of the unit sphere $\Sf^2\subset (P^2)^\perp$ in case $(i)$ and of the hyperbolic plane $\mathbb{H}^2\subset (P^2)^\perp$
in case $(ii)$. Accordingly, the hypersurface is said to be of spherical or hyperbolic type.

    If $P^2$ is degenerate, then $f$ is a rotation hypersurface of parabolic type  parameterized by 
$$f(s,u)=(\gamma_1(s), \gamma_1(s)u_1, \gamma_1(s)u_2, \gamma_4(s)-\frac{1}{2}\gamma_{1}(s)(u_1^2+u_2^2),\gamma_{5}(s)),$$
 with respect to a pseudo--orthonormal basis $\{e_1, \ldots, e_5\}$ of $\R^5_{s+\epsilon_0}$   such that $\<e_1, e_1\>=0=\<e_4, e_4\>$, $\<e_1, e_4\>=1$,  $\<e_2, e_2\>=1=\<e_3, e_3\>$ and $\<e_5, e_5\>=-2(s+ \epsilon_0)+3$, where  $\gamma(s)=(\gamma_1(s), \gamma_4(s), \gamma_5(s))$ is a unit--speed curve in $\Q_s^2(c)\subset P^3=\spa\{e_1, e_4, e_5\}$. 
 
 In each case, one can compute the principal curvatures of $f$ as in \cite{dcd2} and check that the first equation in (\ref{eq:pcurv})
 is satisfied if and only if  
$\gamma_1''+\tilde c\gamma_1=0$, that is,  $\gamma$ is a $\tilde c$--helix in $\Q_s^2(c)\subset \R^3_{s+\epsilon_0}$.  
%

Under the remaining possibilities for $c, \tilde c, s$ and $\tilde s$, either the same conclusions hold or the bilinear form $\beta_x$ in the proof of Theorem \ref{ngeq4} is everywhere degenerate, in which case there exist normal vector fields $\zeta\in \Gamma(N_{\hat f}M)$ and $\tilde N\in \Gamma(N_{\tilde f}M)$ satisfying $\<\zeta, \zeta\>=\tilde \epsilon=\<\tilde N, \tilde N\>$ and $A_{\zeta}^{\hat f}=A_{\tilde N}^{\tilde f}$, and we obtain as before that $f$ and $\tilde f$ are locally given on an open dense subset  as described in Theorem \ref{prob2n4}.

Finally, if one of the principal curvatures of $f$ is zero, then the preceding argument applies with the roles of $f$ and $\tilde f$ interchanged. \qed

\section{Proof of Theorem \ref{main2}}

\noindent \emph{Proof of Theorem \ref{main2}:} Let $(v,V)$ be the pair associated to $f$. Define
\be\label{tildeVi} \tilde V_j=(-1)^{j+1}\delta_j(v_iV_k-v_kV_i), \,\,\,1\leq i\neq j\neq k\leq 3, \,\,\,i<k.\ee
Then $\tilde V=(\tilde V_1, \tilde V_2, \tilde V_3)$ is the unique vector in $\R^3$, up to sign,  such that $(v, |C|^{-1/2}V, |C|^{-1/2}\tilde V)$ is an  orthonormal basis 
of $\R^3$ with respect to the inner product
\be\label{eq:innerproduct} 
\<(x_1, x_2, x_3), (y_1, y_2, y_3)\>=\sum_{i=1}^3\delta_ix_iy_i.
\ee
Therefore, the matrix $D=(v, |C|^{-1/2}V, |C|^{-1/2}\tilde V)$ satisfies $D\delta D^t=\delta$, where $\delta=\mbox{diag}(\hat \epsilon, C/|C|, -\hat\epsilon C/|C|)$.
 It follows that 
$$\hat \epsilon v_iv_j+ C/|C|^2V_iV_j-\hat\epsilon C/|C|^2\tilde V_i\tilde V_j=0, \,\,\,1\leq i\neq j\leq 3.$$
Multiplying by $\epsilon C$ and using that $\hat \epsilon \epsilon=\tilde \epsilon$ and $\hat \epsilon \epsilon C=\hat \epsilon\epsilon\tilde \epsilon(c-\tilde c)=c-\tilde c$ we obtain
$$(c-\tilde c)v_iv_j+\epsilon V_iV_j-\tilde \epsilon \tilde V_i\tilde V_j=0,$$
or equivalently,
$$cv_iv_j+\epsilon V_iV_j=\tilde c v_iv_j+\tilde \epsilon \tilde V_i\tilde V_j.$$
Substituting the preceding equation into $(v)$ yields
$$\frac{\partial h_{ij}}{\partial u_i}+\frac{\partial h_{ji}}{\partial u_j}+h_{ki}h_{kj}+\tilde \epsilon \tilde V_{i} \tilde V_{j} + \tilde cv_iv_j=0.$$
%
On the other hand, differentiating (\ref{tildeVi}) and using equations $(i)$--$(iv)$ yields
$$\frac{\partial \tilde V_j}{\partial u_i}=h_{ij}\tilde V_i, \,\,\,1\leq i\neq j\leq 3.$$
It follows from Proposition \ref{fund} that there exists a hypersurface  $\tilde{f}\colon M^3 \to \Q_{\tilde s}^4(\tilde{c})$ whose first and second fundamental forms are 
$$
I=\sum_{i=1}^3v_i^2du_i^2 \,\,\,\,\mbox{and}\,\,\,\,\,II=\sum_{i=1}^3\tilde V_iv_idu_i^2,
$$
thus $M^3$ admits an isometric immersion into $\Q_{\tilde s}^{4}(\tilde c)$.

Conversely, let $f\colon M^3 \to \Q_s^4(c)$  be a hypersurface for which there exists an  isometric immersion  $\tilde{f}\colon M^3 \to \Q_{\tilde s}^4(\tilde{c})$.
By Lemma \ref{base-comum}, there exists an  orthonormal frame $\{e_1,e_2,e_{3}\}$ of $M^{3}$ of principal directions of both $f$ and $\tilde f$. 
Let  $\lambda_1, \lambda_2, \lambda_3$ and $\mu_1, \mu_2, \mu_3$ be the principal curvatures of  $f$ and $\tilde f$ correspondent to  $e_1, e_2$ and $e_3$, respectively.
Assume that $\lambda_1< \lambda_2< \lambda_3$, and that the unit normal vector field to $f$ has been chosen so that $\lambda_1<0$.
The Gauss equations for $f$ and $\tilde f$ yield
$$c+\epsilon \lambda_i\lambda_j=\tilde c+\tilde \epsilon \mu_i\mu_j, \,\,\,1\leq i\neq j\leq 3.$$
Thus
\be\label{muimuj}
\mu_i\mu_j=C+\hat \epsilon \lambda_i\lambda_j,\,\,\,C=\tilde \epsilon(c-\tilde c), \,\,\,\,1\leq i\neq j\leq 3.
\ee
It follows that
\be\label{muj}
\mu_j^2=\frac{(C+\hat\epsilon\lambda_j\lambda_i)(C+\hat\epsilon\lambda_j\lambda_k)}{C+\hat\epsilon\lambda_i\lambda_k}, \,\,\,1\leq j\neq i\neq k\neq j\leq 3. 
\ee

The  Codazzi equations for $f$ and $\tilde f$ are, respectively,
\begin{eqnarray}
 e_i(\lambda_j) &=& (\lambda_i - \lambda_j)\<\nabla_{e_j}e_i,e_j\>,\ i \neq j, \label{ch1}\\
  (\lambda_j-\lambda_k)\<\nabla_{e_i}e_j,e_k \> &=& (\lambda_i - \lambda_k)\<\nabla_{e_j}e_i,e_k \>,\ \ i \neq j \neq k. \label{ch3}
\end{eqnarray}
and
\begin{eqnarray}
 e_i(\mu_j) &=& (\mu_i - \mu_j)\<\nabla_{e_j}e_i,e_j\>,\ i \neq j,\label{ch11}\\
 (\mu_j-\mu_k)\<\nabla_{e_i}e_j,e_k \> &=& (\mu_i - \mu_k)\<\nabla_{e_j}e_i,e_k \>,\ \ i \neq j \neq k. \label{ch13}
\end{eqnarray}

   Multiplying  (\ref{ch13}) by  $\mu_j$ and using (\ref{muj}) and (\ref{ch3}) we obtain
$$
 \hat \epsilon C\frac{(\lambda_i-\lambda_j)(\lambda_j-\lambda_k)}{C+\hat \epsilon\lambda_i\lambda_k}\<\nabla_{e_i}e_j,e_k \> = 0,\ i \neq j\neq k.
$$
Since the principal curvatures  $\lambda_1, \lambda_2$ and $\lambda_3$ are distinct, it follows that   
\be\label{ijkdistinct}
\<\nabla_{e_i}e_j,e_k \> = 0,\,\,\, 1\leq i \neq j\neq k\neq i\leq 3.
\ee

Computing $2\mu_j e_i(\mu_j)$,  first by differentiating (\ref{muj})  and then by multiplying (\ref{ch11}) by $2\mu_j$, and using (\ref{ch1}), (\ref{muimuj}) and (\ref{muj}), we obtain
\be\label{edp}
 \begin{array}{l} (C+\hat \epsilon\lambda_j\lambda_k)(\lambda_k-\lambda_j)e_i(\lambda_i) + (C+\hat \epsilon\lambda_i\lambda_k)(\lambda_k-\lambda_i)e_i(\lambda_j)\vspace{1ex}\\ \hspace*{20ex}+ (C+\hat \epsilon\lambda_i\lambda_j)(\lambda_i-\lambda_j)e_i(\lambda_k)=0. \end{array}
\ee

Now let $\{\omega_1, \omega_2, \omega_3\}$ be the dual frame of $\{e_1, e_2, e_3\}$, and define the one-forms $\gamma_j$, $1\leq j\leq 3$, by
$$\gamma_j=\sqrt{\delta_j\frac{(\lambda_j-\lambda_i)(\lambda_j-\lambda_k)}{C+\hat\epsilon\lambda_i\lambda_k}}\omega_j, \,\,\,1\leq j\neq i\neq k\neq j\leq 3,$$
where $\delta_j=y_j/|y_j|$ for $y_j=\frac{(\lambda_j-\lambda_i)(\lambda_j-\lambda_k)}{C+\hat\epsilon\lambda_i\lambda_k}$. 

By (\ref{muj}), either all the three numbers $C+\hat\epsilon\lambda_j\lambda_i$, $C+\hat\epsilon\lambda_j\lambda_k$ and $C+\hat\epsilon\lambda_i\lambda_k$ are positive
or two of them are negative and the remaining one is positive. 
 Hence there are four possible cases:
 \begin{itemize}
 \item[(I)] $C+\hat\epsilon \lambda_i\lambda_j>0$, $1\leq i\neq j\leq 3$.
 \item[(II)] $C+\hat\epsilon\lambda_1\lambda_2<0$, $C+\hat\epsilon\lambda_1\lambda_3<0$ and $C+\hat\epsilon\lambda_2\lambda_3>0$.
\item[(III)]   $C+\hat\epsilon\lambda_1\lambda_2>0$, $C+\hat\epsilon\lambda_1\lambda_3<0$ and $C+\hat\epsilon\lambda_2\lambda_3<0$.
 \item[(IV)]  $C+\hat\epsilon\lambda_1\lambda_2<0$, $C+\hat\epsilon\lambda_1\lambda_3>0$ and $C+\hat\epsilon\lambda_2\lambda_3<0$.
 \end{itemize}

Notice that $(\delta_1, \delta_2, \delta_3)$ equals $(1,-1, 1)$ in case $(I)$,  $(1,1, -1)$ in case $(II)$,  $(-1,1, 1)$ in case $(III)$ and $(-1,-1, -1)$ in case $(IV)$. It is easily checked that one must have $\hat \epsilon =-1$ and $C<0$ in case $(IV)$, whereas
in the remaining cases either $\hat \epsilon =1$ or $\hat \epsilon =-1$ and $C>0$. Therefore, $(\delta_1, \delta_2, \delta_3)=
(-1,-1, -1)$ if $\hat \epsilon =-1$ and $C<0$, and in the remaining cases we may assume that 
$(\delta_1, \delta_2, \delta_3)=(1,-1, 1)$  after possibly reordering the coordinates.


We claim that (\ref{edp}) are precisely the conditions for the 
one-forms $\gamma_j$, $1\leq j\leq 3$, to be closed. To prove this, set $x_j=\sqrt{\delta_jy_j}$, $1\leq j\leq 3$, so that $\gamma_j=x_j\omega_j$. It follows from (\ref{ijkdistinct}) that 
$$d\gamma_j(e_i,e_k)=e_i\gamma_j(e_k)-e_k\gamma_j(e_i)-\gamma_j([e_i,e_k])=0.$$
On the other hand, using (\ref{ch1}) we obtain
\begin{eqnarray*}
d\gamma_j(e_i,e_j)&=&e_i\gamma_j(e_j)-e_j\gamma_j(e_i)-\gamma_j([e_i,e_j])\vspace{1ex}\\ 
&=& e_i(x_j)+x_j\<\nabla_{e_j}e_i,e_j\> \vspace{1ex}\\
&=& e_i(x_j)+x_j\dfrac{e_i(\lambda_j)}{\lambda_i-\lambda_j},
\end{eqnarray*}
hence closedness of $\gamma_j$ is equivalent to
\begin{eqnarray}\label{closed}
e_i(x_j)=\dfrac{x_j}{\lambda_j-\lambda_i}e_i(\lambda_j),\,\,\,\  1\leq i\neq j \leq 3.
\end{eqnarray}
We have
$$e_i(x_j)=e_i((\delta_jy_j)^{1/2})=\frac{1}{2}(\delta_jy_j)^{-1/2}\delta_je_i(y_j)=\frac{\delta_j}{2x_j}e_i(y_j),$$
thus (\ref{closed}) is equivalent to
$$\frac{2x_j^2}{\lambda_j-\lambda_i}e_i(\lambda_j)=\delta_je_i(y_j),$$
or yet, to
$$2(\lambda_j-\lambda_k)e_i(\lambda_j)=e_i(y_j)(C+\hat \epsilon\lambda_i\lambda_k).$$
The preceding equation is in turn equivalent to
$$\begin{array}{l}2(\lambda_j-\lambda_k)(C+\hat\epsilon\lambda_i\lambda_k)e_i(\lambda_j) = (e_i(\lambda_j)-e_i(\lambda_i)(\lambda_j-\lambda_k)(C+\hat \epsilon\lambda_i\lambda_k)  \vspace{1ex}\\\hspace*{25ex}+(\lambda_j-\lambda_i)(e_i(\lambda_j)-e_i(\lambda_k))(C+\hat\epsilon\lambda_i\lambda_k)\vspace{1ex}\\\hspace*{25ex}-(\lambda_j-\lambda_i)(\lambda_j-\lambda_k)(\hat\epsilon(e_i(\lambda_i)\lambda_k+\lambda_ie_i(\lambda_k)),\end{array}
$$
which is the same as  (\ref{edp}). 

Therefore, each point $x\in M^3$ has an open neigborhood $V$ where one can find functions $u_j\in C^{\infty}(V)$, $1\leq j\leq 3$, such that $du_j=\gamma_j$, and we can choose
$V$ so small that $\Phi=(u_1, u_2, u_3)$ is a diffeomorphism of $V$ onto an open subset $U\subset \R^3$, that is, $(u_1, u_2, u_3)$ are local coordinates on $V$. From $\delta_{ij}=du_j(\d/\d u_i)=x_j\omega_j(\d/\d u_i)$ it follows that 
$\d/\d u_i=v_ie_i$, with $v_i=x_i^{-1}$.

Now notice that
$$\sum_{j=1}^3\delta_jv_j^2={\sum_{i,k\neq j=1}^3}\dfrac{C+\hat\epsilon\lambda_i\lambda_k}{(\lambda_j-\lambda_i)(\lambda_j-\lambda_k)}=\hat\epsilon,
$$
$$\sum_{j=1}^3\delta_jv_jV_j=\sum_{j=1}^3\delta_j\lambda_jv^2_j={\sum_{i,k\neq j=1}^3}\lambda_j\dfrac{C+\hat\epsilon\lambda_i\lambda_k}{(\lambda_j-\lambda_i)(\lambda_j-\lambda_k)}=0
$$
and 
$$\sum_{j=1}^3\delta_jV_j^2=\sum_{j=1}^3\delta_j\lambda^2_jv^2_j={\sum_{i,k\neq j=1}^3}\lambda^2_j\dfrac{C+\hat\epsilon\lambda_i\lambda_k}{(\lambda_j-\lambda_i)(\lambda_j-\lambda_k)}=C.
$$
It follows that the  pair $(v, V)$ satisfies (\ref{sphol}). \qed

\section{Proof of Theorem \ref{cfhyp}}

  Before starting the proof of Theorem \ref{cfhyp},  recall that the  \emph{Weyl tensor} of a Riemannian manifold $M^n$ is defined by
$$\begin{array}{l} \<C(X,Y)Z, W\>=\<R(X,Y)Z, W\>- L(X,W)\<Y, Z\>- L(Y, Z)\<X,W\>\vspace{1ex}\\
\hspace*{18ex}+L(X, Z)\<Y,W\>+L(Y,W)\<X, Z\>
\end{array}
$$
for all $X, Y, Z, W\in \mathfrak{X}(M)$, where $L$ is the \emph{Schouten tensor} of $M^n$, which is given in terms of the Ricci tensor and the scalar curvature $s$  by
$$L(X, Y)=\frac{1}{n-2}(\Ric(X, Y)-\frac{1}{2}ns\<X, Y\>).$$
It is well-known that, if $n\geq 4$, then the  vanishing of the Weyl tensor is a necessary and sufficient condition for 
 $M^n$ to be conformally flat. \vspace{1ex}\\
\noindent \emph{Proof of Theorem \ref{cfhyp}:} Let $f\colon M^{n} \to \Q_s^{n+1}(c)$ be a conformally flat hypersurface of dimension $n\geq 4$. 
For a fixed point $x\in M^n$, choose a unit normal vector $N\in N_x^fM$ and let $A=A_N\colon T_xM\to T_xM$ 
be the shape operator of $f$ with respect to $N$. Let $W^3$ be a vector space endowed with the Lorentzian inner product $\<\!\<\,,\,\>\!\>$ given by
$$\<\!\<(a,b, c), (a',b',c')\>\!\>=\epsilon(-aa'+bb'+\epsilon cc').$$
Define a bilinear form  $\beta\colon T_xM \times T_xM \to W^3$ by 
$$
\beta(X,Y) =
(L(X,Y)+\dfrac{1}{2}(1-c)\<X,Y\>, L(X,Y)-\dfrac{1}{2}(1+c)\<X,Y\>, \<AX, Y\>).
$$
Note that $\beta(X,X)\neq 0$ for all $X\neq 0$. Moreover, 
$$\begin{array}{l}
\<\!\<\beta(X,Y),\beta(Z,W)\>\!\> - \<\!\<\beta(X,W),\beta(Z,Y)\>\!\>=-L(X,Y)\<Z,W\>\vspace{1ex}\\
\hspace*{2.9ex} - L(Z,W)\<X,Y\> + L(X,W)\<Z,Y\> + L(Z,Y)\<X,W\> + c\<(X \wedge Z)W,Y\> \vspace{1ex}\\
\hspace*{2.9ex}+ \epsilon\<(AX\wedge AZ)W, Y\>  =\<C(X,Z)W, Y\>=0.
\end{array}
$$
Thus $\beta$ is flat with respect to $\<\!\<\,,\,\>\!\>$. 
We claim that $S(\beta)$ must be degenerate. 
Otherwise, we would have 
$$0=\dim \ker \beta\geq n-\dim S(\beta)>0,$$
a contradiction. Now let $\zeta\in S(\beta)\cap S(\beta)^\perp$ and choose a pseudo-orthonormal basis
$\zeta, \eta, \xi$ of $W^3$ with $\<\!\<\zeta, \zeta\>\!\>=0=\<\!\<\eta,\eta\>\!\>$, $\<\!\<\zeta, \eta\>\!\>=1=\<\!\<\xi, \xi\>\!\>$ and
$\<\!\<\xi, \zeta\>\!\>=0=\<\!\<\xi, \eta\>\!\>$. 
Then 
$$\beta=\phi \zeta+ \psi \xi,$$
where $\phi=\<\!\<\beta, \eta\>\!\>$ and $\psi=\<\!\<\beta, \xi\>\!\>$. Flatness of $\beta$ implies that 
$\dim \ker\psi=n-1$. We claim that $\ker \psi$ is an eigenspace of $A$. Given $Z\in \ker \psi$ we have
\be\label{eq:betazx}\beta(Z, X)=\phi(Z,X)\zeta\ee
for all $X\in T_xM$. Let $\{e_1=(1,0,0), e_2=(0,1,0), e_3=(0,0,1)\}$ be the canonical basis of $W$ and write
$\zeta=\sum_{j=1}^3a_je_j$. Then (\ref{eq:betazx}) gives
$$L(Z,X)+\frac{1}{2}(1-c)\<Z, X\>=a_1\phi(Z,X)$$
and
 $$L(Z,X)-\frac{1}{2}(1+c)\<Z, X\>=a_2\phi(Z,X).$$
 Subtracting the second of the preceding equations from the first yields
 $$\<Z, X\>=(a_1-a_2)\phi(Z,X),$$
 which implies that $a_1-a_2\neq 0$ and
 $$\phi(Z,X)=\frac{1}{a_1-a_2}\<Z, X\>.$$
 Moreover, we also obtain from (\ref{eq:betazx}) that
 $$\<AZ, X\>=a_3 \phi(Z,X)=\frac{a_3}{a_1-a_2}\<Z, X\>,$$
 which proves our claim.\vspace{2ex}\qed

\section{Proof of Theorem \ref{main3}}

 First  recall that 
a necessary and sufficient condition for  a three-dimensional Riemannian manifold $M^3$ 
to be conformally flat is that its Schouten tensor $L$ be a \emph{Codazzy tensor}, that is,
$$(\nabla_X L)(Y,Z)=  (\nabla_Y L)(X,Z)$$
for all $X, Y, Z\in \mathfrak{X}(M)$, where 
$$ (\nabla_X L)(Y,Z)=X(L(Y,Z))-L(\nabla_X Y, Z)-L(Y, \nabla_XZ).$$

\noindent \emph{Proof of Theorem \ref{main3}:} Let $f\colon\,M^3\to\Q_s^{4}(c)$ 
be a  holonomic hypersurface whose associated pair $(v, V)$ satisfies 
(\ref{eq:cflatcond}).  Then $v=(v_1, v_2, v_3)$ is a null vector with respect to the Lorentzian inner product $\<\, ,\,\>$ given by 
(\ref{eq:innerproduct}), with $(\delta_1, \delta_2, \delta_3)=(1, -1, 1)$,  and $V=(V_1, V_2, V_3)$ is a unit space-like vector orthogonal to $v$. Thus, we may write
$$V=\frac{\rho}{v_2}v+\frac{\lambda}{v_2}(-v_3,0,v_1), \,\,\,\,\lambda=\pm 1,$$
for some $\rho\in C^{\infty}(M)$, which is equivalent to
\be\label{eq:V1V3}
V_1=\frac{1}{v_2}(V_2v_1-\lambda v_3)\,\,\,\,\mbox{and}\,\,\,\,V_3=\frac{1}{v_2}(V_2v_3+\lambda v_1).
\ee
The eigenvalues $\mu_1, \mu_2$ and $\mu_3$ of the  Schouten tensor $L$  are given by
$$2\mu_j=c+\epsilon(\lambda_i\lambda_j+\lambda_k\lambda_j-\lambda_i\lambda_k),\,\,\,1\leq j\leq 3,$$
where $\lambda_j$, $1\leq j\leq 3$, are the principal curvatures of $f$.
Define
\be\label{eq:phij}
\phi_j=v_j(\lambda_i\lambda_j+\lambda_k\lambda_j-\lambda_i\lambda_k),\,\,\,1\leq j\leq 3.
\ee
That  $L$ is a Codazzi tensor is then equivalent to  the equations
\be\label{eq:codazziphi}\frac{\partial \phi_j}{\partial u_i}=h_{ij}\phi_i,\,\,\,1\leq i\neq j\leq 3.\ee
Replacing $\lambda_j=\frac{V_j}{v_j}$ in (\ref{eq:phij}) and using (\ref{eq:V1V3}) we obtain
$$\phi_1=\frac{1}{v_2^2}(-2\lambda V_2v_3+(V_2^2-1)v_1),\,\,\,\,\,\,\,\phi_2=\frac{1}{v_2}(V_2^2+1)$$
and
$$\phi_3=\frac{1}{v_2^2}((V_2^2-1)v_3+2\lambda V_2v_1).$$
It is now a straightforward computation to verify (\ref{eq:codazziphi}) by using equations $(i)$ and $(iv)$ of system (\ref{sistema-hol}) together with equations 
(\ref{eq:int1}) and (\ref{eq:int2}).

Conversely, assume that $f\colon\,M^3\to\Q_s^{4}(c)$  is an isometric immersion with three distinct principal curvatures 
$\lambda_1<\lambda_2<\lambda_3$ of a conformally flat manifold. Let $\{e_1,e_2,e_3\}$ be a correspondent  orthonormal frame of principal directions. Then  $\{e_1,e_2,e_3\}$ also diagonalyzes the Schouten tensor $L$, 
 and the correspondent  eigenvalues are  
\begin{equation}\label{muj2}
2\mu_j = \epsilon(\lambda_i\lambda_j + \lambda_j\lambda_k - \lambda_i\lambda_k) + c, \,\,\,\,1\leq j\leq 3.
\end{equation}

The  Codazzi equations for $f$ and $L$ are, respectively,
\begin{eqnarray}
 e_i(\lambda_j) &=& (\lambda_i - \lambda_j)\<\nabla_{e_j}e_i,e_j\>, \,\,\,\ i \neq j, \label{ch1b}\\
  (\lambda_j-\lambda_k)\<\nabla_{e_i}e_j,e_k \> &=& (\lambda_i - \lambda_k)\<\nabla_{e_j}e_i,e_k \>,\,\,\,\ \ i \neq j \neq k. \label{ch3b}
\end{eqnarray}
and
\begin{eqnarray}
 e_i(\mu_j) &=& (\mu_i - \mu_j)\<\nabla_{e_j}e_i,e_j\>,\,\,\,\ i \neq j,\label{ch11b}\\
 (\mu_j-\mu_k)\<\nabla_{e_i}e_j,e_k \> &=& (\mu_i - \mu_k)\<\nabla_{e_j}e_i,e_k \>,\,\,\,\ \ i \neq j \neq k. \label{ch13b}
\end{eqnarray}

Substituting $(\ref{muj2})$ into $(\ref{ch13b})$, and using  $(\ref{ch3b})$, we obtain
\begin{equation}
(\lambda_i-\lambda_j)(\lambda_j-\lambda_k)\<\nabla_{e_i}e_j,e_k\>=0,\,\,\,\,\,\, i \neq j \neq k.
\end{equation}
Since  $\lambda_1,\lambda_2$ and $\lambda_3$ are pairwise distinct, it follows that
\begin{equation}
\<\nabla_{e_i}e_j,e_k\>=0,\,\,\,\,\,\, 1\leq i \neq j \neq k\neq i \leq 3.
\end{equation}

Differentiating  $(\ref{muj2})$ with respect to  $e_i$, we obtain
\begin{eqnarray}\label{diff}
2e_i(\mu_j)=\epsilon[(\lambda_i+\lambda_k)e_i(\lambda_j) + (\lambda_j-\lambda_k)e_i(\lambda_i) + (\lambda_j-\lambda_i)e_i(\lambda_k)].
\end{eqnarray}
On the other hand, it follows from  $(\ref{ch1})$, $(\ref{ch11b})$ and $(\ref{muj2})$ that
\begin{eqnarray}\label{usoutd}
e_i(\mu_j)=\epsilon\lambda_ke_i(\lambda_j).
\end{eqnarray}
Hence
\begin{equation}\label{edp2}
(\lambda_j - \lambda_k)e_i(\lambda_i)+(\lambda_i - \lambda_k)e_i(\lambda_j)+(\lambda_j - \lambda_i)e_i(\lambda_k)=0. 
\end{equation}

Now let $\{\omega_1, \omega_2, \omega_3\}$ be the dual frame of $\{e_1, e_2, e_3\}$, and define the one-forms $\gamma_j$, $1\leq j\leq 3$, by
\be\label{gammajb} 
\gamma_j=\sqrt{\delta_j(\lambda_j-\lambda_i)(\lambda_j-\lambda_k)}\,\omega_j, \,\,\,\,\,\,1\leq j\neq i\neq k\neq j\leq 3,
\ee
where $(\delta_1, \delta_2, \delta_3)=(1, -1, 1)$. One can check that  (\ref{edp2}) are precisely the conditions for the 
one-forms $\gamma_j$, $1\leq j\leq 3$, to be closed. 

Therefore, each point $x\in M^3$ has an open neigborhood $V$ where one can find functions $u_j\in C^{\infty}(V)$, $1\leq j\leq 3$, such that $du_j=\gamma_j$, and we can choose
$V$ so small that $\Phi=(u_1, u_2, u_3)$ is a diffeomorphism of $V$ onto an open subset $U\subset \R^3$, that is, $(u_1, u_2, u_3)$ are local coordinates on $V$. From $\delta_{ij}=du_j(\partial_{\,i})=x_j\omega_j(\partial_{\,i})$ it follows that 
$\partial_{j}=v_je_j$, $1\leq j\leq 3$, with 
$$v_j=\sqrt{\frac{\delta_j}{(\lambda_j-\lambda_i)(\lambda_j-\lambda_k)}}$$.

Now notice that
$$\sum_{j=1}^3\delta_jv_j^2={\sum_{i,k\neq j=1}^3}\dfrac{1}{(\lambda_j-\lambda_i)(\lambda_j-\lambda_k)}=0,$$
$$\sum_{j=1}^3\delta_jv_jV_j=\sum_{j=1}^3\delta_j\lambda_jv^2_j={\sum_{i,k\neq j=1}^3}\dfrac{\lambda_j}{(\lambda_j-\lambda_i)(\lambda_j-\lambda_k)}=0$$
and
$$\sum_{j=1}^3\delta_jV^2_j=\sum_{j=1}^3\delta_j\lambda^2_jv^2_j={\sum_{i,k\neq j=1}^3}\dfrac{\lambda^2_j}{(\lambda_j-\lambda_i)(\lambda_j-\lambda_k)}=1.$$
It follows that $(v, V)$ satisfies 
(\ref{eq:cflatcond}). \qed

\section{The  Ribaucour transformation}

Two  immersions $f\colon M^n \to \R_s^{n+p}$ and $f'\colon M^n \to \R_s^{n+p}$ are said to be related by a Ribaucour transformation if $|f-f'|\neq 0$ everywhere and there exist a vector bundle isometry  $\mathcal{P}\colon f^*T\R_s^{n+p} \to f'^*T\R_s^{n+p}$, a  tensor $D\in \Gamma(T^*M\otimes TM)$, which is symmetric with respect to the induced metrics, 
 and a nowhere vanishing  $\delta\in \Gamma(f^*T\R_s^{n+p})$ such that 
\begin{itemize}
\item[(a)] $\mathcal{P}(Z) - Z= \<\delta, Z\>(f-f')$ for all $Z \in \Gamma(f^*T\R_s^{n+p})$;
\item[(b)] $\mathcal{P} \circ f_* \circ D = f'_*$.
\end{itemize}

Given an  immersion $f\colon M^n \to \Q_s^{n+p}(c)$,  with $c\neq 0$, let $F=i\circ f\colon M^n \to \R_{s+\epsilon_0}^{n+p+1}$, where $\epsilon_0=0$ or $1$ corresponding to $c>0$ or $c<0$, respectively, and $i\colon \Q_s^{n+p}(c)\to \R_{s+\epsilon_0}^{n+p+1}$ denotes an umbilical inclusion. An immersion $f'\colon M^n \to \Q_s^{n+p}(c)$ is said to be a Ribaucour transform of $f$ with data $({\cal P}, D, \delta)$ if $F'=i\circ f'\colon  M^n \to \R_{s+\epsilon_0}^{n+p+1}$ is a Ribaucour transform of $F$ with data $(\hat{\cal P}, D, \hat\delta)$, where $\hat\delta=\delta-cF$ and 
$\hat{\mathcal{P}}\colon F^*T\R_{s+\epsilon_0}^{n+p+1} \to F'^*T\R_{s+\epsilon_0}^{n+p+1}$  is the extension of ${\cal P}$ such that $\hat{\cal P}(F)=F'$. The next result was proved in \cite{dt3}.

\begin{theorem}\label{pr:rcsc}
Let $f\colon\,M^n\to\Q_s^{n+p}(c)$
be an isometric immersion of a simply connected Riemannian manifold and let
$f'\colon\, M^n\to\Q_s^{n+p}(c)$ be a Ribaucour transform of $f$ with data
$({\cal P},D,\delta)$. Then there exist $\varphi\in C^{\infty}(M)$ and $\hat\beta\in \Gamma(N_fM)$  satisfying 
\be\label{vabeta}
\alpha_f(\nabla\varphi, X)+\nabla_X^\perp \hat\beta=0\,\,\,\mbox{for all}\,\,\,X\in TM
\ee
such~that $F'=i\circ f'$ and $F=i\circ f$ are related by 
\be\label{eq:rb2}
F'= F - 2\nu\varphi{\cal G},
\ee
where ${\cal G}=F_*\nabla\varphi + i_*\hat\beta +c\varphi F$ and
$\nu=\<{\cal G},{\cal G}\>^{-1}$.
Moreover,
\be\label{eq:pdo2}
\hat{\cal P}=I - 2\nu{\cal G}{\cal G}^*,\;\;\;
D=I - 2\nu\varphi\Phi
\;\;\;\mbox{and}\;\;\;\hat{\delta}
=-\varphi^{-1}{\cal G},
\ee
where $\Phi=
\hess\varphi +c\varphi I -A^f_{\hat\beta}$. Conversely, given $\varphi\in C^{\infty}(M)$ and $\hat\beta\in \Gamma(N_fM)$  satisfying (\ref{vabeta}) such that $\varphi\nu\neq 0$
everywhere, let $U\subset M^n$
be an open subset where the tensor $D$ given by {\em (\ref{eq:pdo2})} is
invertible, and let
$F'\colon\,U\to\R_{s+\epsilon_o}^{n+p+1}$ be defined by
{\em (\ref{eq:rb2})}. Then $F'=i\circ f'$, where
$f'$ is a Ribaucour transform of $f$. Moreover, the second fundamental forms of $f$ and $f'$ are related by
\be\label{sffsF}
\tilde A^{f'}_{{\cal P}\xi}=D^{-1}(A^f_{\xi}+2\nu\<\hat\beta,\xi\>\Phi), \,\,\,\mbox{for any}\,\,\,\xi\in \Gamma(N_fM).
\ee
\end{theorem}

We now derive from Theorem \ref{pr:rcsc} a Ribaucour transformation for holonomic hypersurfaces,
in a form that is slightly different from  the one in \cite{dt2}. For that we need the following.

\begin{proposition}\label{linearsystem}
Let $f\colon\,M^n\to\Q_s^{n+1}(c)$ 
be a holonomic hypersurface with associated pair $(v,V)$.   Then, the linear system of PDE's
\begin{eqnarray}\label{sist-Rib}
\left\{\begin{array}{l}
(i) \dfrac{\partial \varphi}{\partial u_i} = v_i\gamma_i,\,\,\,\,\,\,\,
(ii) \dfrac{\partial \gamma_j}{\partial u_i} = h_{ji}\gamma_i,\ \ \,\,\ \ i\neq j,\vspace{1ex}\\
(iii) \dfrac{\partial \gamma_i}{\partial u_i} = (v_i - v_i')\psi - \sum_{j \neq i}h_{ji}\gamma_j + \beta V_i - c \varphi v_i, \vspace{1ex}\\
(iv) \epsilon\dfrac{\partial \beta}{\partial u_i} = -V_{i}\gamma_i,\,\,\,\epsilon=-2s+1,\vspace{1ex}\\
(v) \dfrac{\partial \log \psi}{\partial u_i} = -\dfrac{\gamma_i v_i'}{\varphi},\,\,\,\,\,\,
(vi) \dfrac{\partial v_i'}{\partial u_j} = h_{ji}'v_j',\ \ \,\,\ \ i\neq j,
\end{array}\right.
\end{eqnarray}
with $h_{ij}$ and $h'_{ij}$ given by (\ref{eq:hij'}), is completely integrable and has the first integral 
\begin{equation}\label{int-primeira}
\sum_i {\gamma_i}^2+\epsilon\beta^2 + c\varphi^2 - 2\varphi \psi = K\in \R.
\end{equation}
\end{proposition}
\proof A straightforward computation.\qed

\begin{theorem}\label{teo-Rib-bom}
Let $f\colon\,M^n\to\Q_s^{n+1}(c)$ 
be a holonomic hypersurface with associated pair $(v,V)$. If $f'\colon\,M^n\to\Q_s^{n+1}(c)$  is a Ribaucour transform of $f$, then there exists a solution 
$(\gamma,v',\varphi,\psi,\beta)$ of (\ref{sist-Rib}) satisfying 
\begin{equation}\label{int-primeirab}
\sum_i {\gamma_i}^2+\epsilon\beta^2 + c\varphi^2 - 2\varphi \psi = 0
\end{equation}
 such that  $F'=i\circ f'$ and $F=i\circ f$ are related by  
\begin{equation}\label{dt2:10}
{F'} = F - \dfrac{1}{\psi}\left(\sum_i \gamma_i F_*e_i + \beta i_*\xi + c\varphi F \right),
\end{equation}
where $\xi$ is a  unit normal vector field to $f$ and $e_i=v_i^{-1}\partial_{\,i}$, $1\leq i\leq n$. 

Conversely, given a solution  $(\gamma,v',\varphi,\psi,\beta)$ of (\ref{sist-Rib}) satisfying (\ref{int-primeirab}) on an open subset $U\subset M^n$ where
$v_i'$ is positive  for $1\leq i\leq n$, then  $F'$ defined by (\ref{dt2:10}) is an immersion such that $F'= i \circ f'$, where 
$f'$ is a Ribaucour transform of $f$ whose associated pair is $(v', V')$, with
\be\label{asstriple}
V_i' = V_i + (v_i - v_i')\dfrac{\epsilon\beta}{\varphi}, \,\,\,1\leq i\leq n.
\end{equation}
\end{theorem}
\proof Let $f'\colon\,M^n\to\Q_s^{n+1}(c)$  be  a Ribaucour transform of $f$. By Theorem~\ref{pr:rcsc},
 there exist $\varphi\in C^{\infty}(M)$ and $\hat\beta\in \Gamma(N_fM)$ satisfying (\ref{vabeta}) 
such~that $F'$ is given by (\ref{eq:rb2}), 
where ${\cal G}=F_*\nabla\varphi + i_*\hat\beta +c\varphi F$ and
$\nu=\<{\cal G},{\cal G}\>^{-1}$. 

Write $\nabla\varphi=\sum_{i=1}^n\gamma_ie_i$, where $ \gamma_i\in C^{\infty}(M)$, $1\leq i\leq n$. 
Since $\partial_{\,i}=v_ie_i$, $1\leq i\leq n$, this is equivalent to equation $(i)$
of system (\ref{sist-Rib}). Now write $\hat \beta=\beta\xi$, where $\beta\in C^{\infty}(M)$.
Then (\ref{vabeta}) can be written as 
\be\label{vabetab}
A\nabla\varphi=-\epsilon \nabla\beta,
\ee
which is  equivalent, by taking inner products of both sides with  $\partial_i$,
to equation $(iv)$
of system (\ref{sist-Rib}). On the other hand, equation (\ref{vabeta}) implies that 
$$ {\cal G}_*=F_*\Phi,$$
where $\Phi=
\hess\varphi +c\varphi I -A^f_{\hat\beta}$. 
Therefore $\Phi$ is a Codazzi tensor that satisfies
$$\alpha_f(\Phi X, Y)=\alpha_f(X, \Phi Y)$$
for all $X, Y\in TM$, that is, $\Phi$ has $\{e_1, \ldots, e_n\}$ as a diagonalyzing frame.
Since 
\begin{equation} \label{phi-bom}
\Phi \partial_{\,i} = \left(\dfrac{\partial \gamma_i}{\partial u_i} + \sum_{j \neq i}h_{ji}\gamma_j - \beta V_i + c v_i \varphi \right)e_i + \sum_{j \neq i}\left(\dfrac{\partial \gamma_j}{\partial u_i} - h_{ji}\gamma_i \right)e_j,
\end{equation}
equation $(ii)$ of system (\ref{sist-Rib}) follows.

Now define $\psi\in C^{\infty}(M)$ by 
$$2\varphi\psi=\<{\cal G},{\cal G}\>=\sum_i {\gamma_i}^2+\epsilon\beta^2 + c\varphi^2.$$
Differentiating both sides with respect to $u_i$ and using equations $(i)$, $(ii)$  and $(iv)$ of (\ref{sist-Rib}) yields
\begin{eqnarray}\label{derivada-int1}
 \dfrac{\partial \gamma_i}{\partial u_i} + \sum_{j \neq i}h_{ji}\gamma_j - \beta V_i + c v_i \varphi = v_i \psi + \dfrac{\varphi}{\gamma_i}\dfrac{\partial \psi}{\partial u_i}.
\end{eqnarray}
Defining $v_i'$ by  $(v)$, then $(iii)$ follows from (\ref{derivada-int1}).

Finally, from $\displaystyle{\dfrac{\partial^2 \gamma_i}{\partial u_i \partial u_j} = \dfrac{\partial^2 \gamma_i}{\partial u_j \partial u_i}}$ we obtain
\begin{eqnarray*}
&&\dfrac{\partial}{\partial u_i}\left(h_{ij}\gamma_j\right) = \dfrac{\partial}{\partial u_j}\left((v_i - v_i')\psi - \sum_{k \neq i}h_{ki}\gamma_k + \beta V_i - c\varphi v_i \right),
\end{eqnarray*}
thus
$$
\begin{array}{l}
\dfrac{\partial h_{ij}}{\partial u_i}\gamma_j + h_{ij}\dfrac{\partial \gamma_j}{\partial u_i} = \left(\dfrac{\partial v_i}{\partial u_j} - \dfrac{\partial v_i'}{\partial u_j}\right)\psi + (v_i - v_i')\dfrac{\partial \psi}{\partial u_j} - \dfrac{\partial h_{ji}}{\partial u_j}\gamma_j - h_{ji}\dfrac{\partial \gamma_j}{\partial u_j} \vspace{1ex}\\\hspace*{14ex}- \dfrac{\partial h_{ki}}{\partial u_j}\gamma_k  - h_{ki}\dfrac{\partial \gamma_k}{\partial u_j} + \dfrac{\partial \beta}{\partial u_j}V_i + \beta \dfrac{\partial V_i}{\partial u_j} - c\dfrac{\partial \varphi}{\partial u_j}v_i - c\varphi \dfrac{\partial v_i}{\partial u_j}.
\end{array}
$$
It follows that
$$
\begin{array}{l}
\left(\dfrac{\partial h_{ij}}{\partial u_i} + \dfrac{\partial h_{ji}}{\partial u_j}\right)\gamma_j  + h_{ij}h_{ji}\gamma_i =
\psi h_{ji}v_j - \dfrac{\partial v_i'}{\partial u_j}\psi - (v_i - v_i')\dfrac{\gamma_j \psi}{\varphi}v_j' \vspace{1ex}\\ - h_{ji}(v_j - v_j')\psi  + h_{ji}h_{ij}\gamma_i + h_{ji}h_{kj}\gamma_k - \beta h_{ji}V_j + c \varphi h_{ji}v_j - h_{kj}h_{ji}\gamma_k \vspace{1ex}\\ - h_{ki}h_{kj}\gamma_j - V_iV_j\gamma_j  + \beta h_{ji}V_j - cv_iv_j\gamma_j - c\varphi h_{ji}v_j,
\end{array}
$$
which yields equation $(vi)$ of (\ref{sist-Rib}). 

Conversely, let $F'$ be given by (\ref{dt2:10}) in terms of  a solution  $(\gamma,v',\varphi,\psi,\beta)$ of (\ref{sist-Rib}) satisfying (\ref{int-primeirab}) on an open subset $U\subset M^n$ where
$v_i'$ is nowhere vanishing  for $1\leq i\leq n$. We have $\nabla\varphi=\sum_{i=1}^n\gamma_ie_i$ by equation $(i)$ of (\ref{sist-Rib}). Defining $\hat \beta\in \Gamma(N_fM)$ by $\hat \beta=\beta\xi$, we can write $F'$ as in (\ref{eq:rb2}), with ${\cal G}=F_*\nabla\varphi + i_*\hat\beta +c\varphi F$ and
$\nu=\<{\cal G},{\cal G}\>^{-1}$. In view of $(iv)$, equation (\ref{vabetab}) is satisfied, and hence so  is (\ref{vabeta}).
Thus ${\cal G}_*=F_*\circ \Phi$, where  $\Phi=
\hess\varphi +c\varphi I -A^f_{\hat\beta}$. 

It follows from $(ii)$ and (\ref{phi-bom}) that $\Phi \partial_{\,i}=B_i\partial_i,$
where $$B_i=v_i^{-1}\left(\dfrac{\partial \gamma_i}{\partial u_i} + \sum_{j \neq i}h_{ji}\gamma_j - \beta V_i + c v_i \varphi\right)=v_i^{-1}(v_i-v_i')\psi.$$
Using $(iii)$ and (\ref{int-primeirab}) we obtain 
$$D\partial_{\,i}=(1-2\nu\varphi B_i)\partial_i=(1-2\nu\varphi v_i^{-1}(v_i-v_i')\psi)=\frac{v_i'}{v_i}\partial_{\,i}.
$$
Thus $D$ is invertible wherever $v_i'$ does not vanish  for $1\leq i\leq n$.  It follows from Theorem 
\ref{pr:rcsc} that the map  $F'$ defined by (\ref{dt2:10}) is an immersion on $U$ and that $F'= i \circ f'$,
where $f'$ is a Ribaucour transform of $f$. Moreover, we obtain from (\ref{sffsF}) that $F'$, and hence $f'$, is 
holonomic with $u_1, \ldots, u_n$ as principal coordinates. It also follows from (\ref{sffsF}) that
$$\frac{V_i'}{v_i'}{\partial_{\,i}}= A^{f'}{\partial_i}=\frac{v_i}{v_i'}\left(\frac{V_i}{v_i}+\frac{\epsilon\beta}{\varphi}\frac{v_i-v'_i}{v_i}\right){\partial_{\,i}},$$
which yields (\ref{asstriple}).\qed


\subsection{The Ribaucour transformation for solutions of \\ Problem~$*$ and for conformally flat hypersurfaces.}

We now specialize the Ribaucour transformation to the classes of hypersurfaces $f\colon\,M^3\to\Q_s^{4}(c)$ 
 that are either conformally flat or admit an isometric immersion into  $\Q_{\tilde s}^4(\tilde c)$ with $\tilde c\neq c$.

\begin{proposition}\label{speciallinearsystem}
Let $f\colon\,M^3\to\Q_s^{4}(c)$  be a  holonomic  hypersurface 
  whose associated pair $(v,V)$ satisfies (\ref{sphol}) (respectively, (\ref{eq:cflatcond})). Then, the linear system of PDE's
obtained by adding the equation 
\be\label{vii}
 \delta_i\dfrac{\partial v_i'}{\partial u_i}+ \delta_j h'_{ij}v_j'+ \delta_k h'_{ik}v_k'=0
\ee
to system (\ref{sist-Rib}), 
where $h'_{ij}$ is given by (\ref{eq:hij'}), 
 is completely integrable and has (besides (\ref{int-primeira})) the first integral
\begin{equation}\label{int-primeira2}
\delta_1 v_1'^2 + \delta_2 v_2'^2 + \delta_3 v_3'^2 = K\in \R. 
\end{equation}
Moreover, 
the function 
\be\label{omega}
\Omega=\varphi\sum_{j=1}^3\delta_j v_j'V_j  - \epsilon\beta \left(K -\sum_{j=1}^3\delta_j  v_jv_j'\right)
\ee
satisfies
\be\label{eqomega}
\dfrac{\d \Omega}{\d u_i}=\dfrac{\gamma_i}{\varphi}(v_i+v'_i)\Omega.
\ee
In particular, if initial conditions for $\varphi$ and $\beta$ at $x_0\in M^3$ are chosen so that $\Omega$ vanishes at $x_0$,
 then $\Omega$ vanishes everywhere.
\end{proposition}
\proof  The first two assertions follow from straightforward computations. To prove the last one,  define $\rho=\sum_{i=1}^3 \delta_iv_i'V_i$ and $\Theta=K - \sum_{i=1}^3 \delta_iv_i'v_i$. We have
$$\begin{array}{l}\dfrac{\d \rho}{\d u_i}= \delta_i\dfrac{\d v'_i}{\d u_i}V_i+\delta_iv_i'\dfrac{\d V_i}{\d u_i}+\sum_{j\neq i}\delta_j\dfrac{\d v'_j}{\d u_i}V_j+\sum_{j\neq i}\delta_jv'_j\dfrac{\d V_j }{\d u_i}\vspace{1ex}\\
\hspace*{4.5ex}=\sum_{j\neq i}\delta_j(h_{ij}-h'_{ij}) v'_jV_i-\sum_{j\neq i}\delta_j(h_{ij}-h'_{ij}) V_jv'_i\vspace{1ex}\\
\hspace*{4.5ex}=
\sum_{j\neq i}\delta_j(v'_j -v_j) V_j\dfrac{\gamma_iv'_i}{\varphi}-\sum_{j\neq i}\delta_j(v'_j -v_j) v'_j\dfrac{\gamma_iV_i}{\varphi}\vspace{1ex}\\
\hspace*{4.5ex}=\dfrac{v'_i\gamma_i}{\varphi}\left(\rho-\delta_iv'_iV_i+\delta_iv_iV_i\right)-\dfrac{V_i\gamma_i}{\varphi}\left(\Theta-\delta_i{v'_i}^2 +\delta_iv_iv'_i)\right)\vspace{1ex}\\
\hspace*{4.5ex}= \dfrac{\gamma_i}{\varphi}(v'_i\rho-\Theta V_i)\end{array}
$$
and 
\begin{eqnarray*}\dfrac{\d \Theta}{\d u_i}&=& -\delta_i\dfrac{\d v_i}{\d u_i}v'_i-\delta_iv_i\dfrac{\d v'_i}{\d u_i}-\sum_{j\neq i}\delta_j\dfrac{\d v_j}{\d u_i}v'_j-\sum_{j\neq i}\delta_jv_j\dfrac{\d v'_j }{\d u_i}\\
 &=&\left(\sum_{j\neq i}\delta_jv_j(h_{ij}-h'_{ij})\right) v'_i+\left(\sum_{j\neq i}\delta_j(h'_{ij}-h_{ij})v'_j\right) v_i\\
&=&\left(\sum_{j\neq i}\delta_jv_j(v_j-v_j')\right)\dfrac{\gamma_iv_i'}{\varphi}+\left(\sum_{j\neq i}\delta_j(v'_j-v_j)v'_j\right)\dfrac{v_i\gamma_i}{\varphi}\\
&=& (\Theta-\delta_iv_i^2+\delta_iv_iv'_i))\dfrac{\gamma_iv_i'}{\varphi}+
(\Theta-\delta_i{v'_i}^2+\delta_iv_iv'_i))\dfrac{v_i\gamma_i}{\varphi}\\
&=& \dfrac{\gamma_i}{\varphi}(v_i+v'_i)\Theta.
\end{eqnarray*}
Therefore,
\begin{eqnarray*}
\dfrac{\d \Omega}{\d u_i}&=&\dfrac{\d \varphi}{\d u_i}\rho+\varphi\dfrac{\d \rho}{\d u_i}-\dfrac{\d \epsilon\beta}{\d u_i}\Theta-\epsilon\beta\dfrac{\d \Theta}{\d u_i}\\
&=& v_i\gamma_i\rho+\varphi \dfrac{\gamma_i}{\varphi}(v'_i\rho-\Theta V_i) +V_i\gamma_i\Theta-\epsilon \beta\dfrac{\gamma_i}{\varphi}(v_i+v'_i)\Theta\\
&=& \rho\gamma_i(v_i+v'_i)-\dfrac{\epsilon \beta\gamma_i}{\varphi}(v_i+v'_i)\Theta\\
&=&\dfrac{\gamma_i}{\varphi}(v_i+v'_i)\Omega,
\end{eqnarray*}
which proves (\ref{eqomega}). The last assertion  follows from (\ref{eqomega}) and the lemma below.

\begin{lemma} Let $M^n$ be a connected manifold and let $\Omega\in C^{\infty}(M)$. Assume that there exists a smooth one-form $\omega$ on $M^n$ such that
$d\Omega=\omega \Omega$.
If $\Omega$ vanishes at some point of $M^n$, then it vanishes everywhere.
\end{lemma}
\proof
Given any smooth curve $\gamma\colon I\to M^n$ with $0\in I$, denote $\lambda(s)=\omega(\gamma'(s))$. By the assumption we have 
$$(\Omega\circ \gamma)(t)=(\Omega\circ \gamma)(0)\exp \int_0^t\lambda(s)ds,$$
 and the conclusion follows from the connectedness of $M^n$.\vspace{2ex}\qed

The next result contains Theorem \ref{ribA} in the introduction.

\begin{theorem}\label{thm:ribsol}
Let  $f \colon M^3 \to \Q_s^4(c)$ be a  holonomic  hypersurface whose associated pair $(v,V)$ satisfies (\ref{sphol}) (respectively, (\ref{eq:cflatcond}))
and  $f' \colon M^3 \to \Q_s^4(c)$  a Ribaucour transform of $f$
determined by a solution $(\gamma_1,\gamma_2,\gamma_3,v_1',v_2',v_3',\varphi,\psi,\beta)$ of system (\ref{sist-Rib}). 
If the associated pair $(v', V')$ of  $f'$  also satisfies (\ref{sphol}) (respectively, (\ref{eq:cflatcond})), then 
\be\label{omega3}
\Omega:=\varphi\sum_{j=1}^3\delta_j v_j'V_j  - \epsilon\beta \left(K -\sum_{j=1}^3\delta_j  v_jv_j'\right)=0,
\ee
with $K=\hat\epsilon$ (respectively, $K=0$).
Conversely, let $(\gamma_1,\gamma_2,\gamma_3,v_1',v_2',v_3',\varphi,\psi,\beta)$ be a solution of the linear system of PDE's
obtained by adding  equation (\ref{vii}) to system (\ref{sist-Rib}). If (\ref{int-primeirab}), 
 (\ref{omega3}) and 
\be\label{vi'}
\sum_{i=1}^3\delta_i{v'_i}^2=K,
\ee
where $K=\hat\epsilon$ (respectively, $K=0$), are satisfied at some point of $M^3$, then
(they are satisfied at every point of $M^3$ and) the  pair $(v', V')$ associated to the  Ribaucour transform of $f$
determined by such a solution  also satisfies (\ref{sphol}) (respectively, (\ref{eq:cflatcond})).
\end{theorem}
\proof Let $(v',V')$ be the pair associated to $f'$. Then, using conditions (\ref{sphol}) (respectively, (\ref{eq:cflatcond})),
we obtain
\begin{eqnarray}\label{sums}
\sum_{j=1}^3\delta_j{V_j'}^2-\sum_{j=1}^3\delta_j{V_j}^2&=&\sum_{j=1}^3\delta_j (V_j'-V_j)(V_j'+V_j)\nonumber\\
&=& \frac{\epsilon\beta}{\varphi}\sum_{j=1}^3\delta_j(v_j-v_j')\left(2V_j+\frac{\epsilon\beta}{\varphi}(v_j-v_j')\right)\nonumber\\
&=&\frac{\epsilon\beta}{\varphi}\left(2\sum_{j=1}^3\delta_jV_j(v_j-v_j')+\frac{\epsilon\beta}{\varphi}\sum_{j=1}^3\delta_j(v_j-v_j')^2\right)\nonumber\\
&=&\frac{\epsilon\beta}{\varphi^2}\left(-2\Omega+\epsilon\beta\left(\sum_{j=1}^3\delta_j{v'_j}^2-K\right)\right),
\end{eqnarray}
where $K=\hat \epsilon$ (respectively, $K=0$). If the pair $(v',V')$ associated to  $f'$  satisfies (\ref{sphol}) (respectively, (\ref{eq:cflatcond})), then 
(\ref{vi'}) holds, as well as 
\begin{equation}\label{vi'Vi'}
\sum_{j=1}^3\delta_j v_j'V'_j=0
\end{equation}
and
\begin{equation}\label{cond-desejavel2}
\sum_{j=1}^3\delta_j {V'_j}^2=C,
\end{equation}
where $C=\tilde \epsilon(c-\tilde c)$ (respectively, $C=1$). It follows from (\ref{sums}) that 
(\ref{omega3}) holds. 
 
 Conversely, let $(\gamma_1,\gamma_2,\gamma_3,v_1',v_2',v_3',\varphi,\psi,\beta)$ be a solution of the linear system of PDE's
obtained by adding  equation (\ref{vii}) to system (\ref{sist-Rib}). If (\ref{int-primeirab}), (\ref{vi'}) and (\ref{omega3}) are satisfied at some point of $M^n$, then
they are satisfied at every point of $M^n$ by Proposition \ref{speciallinearsystem}. Then, equations (\ref{vi'}), (\ref{omega3}) and (\ref{sums})
 imply  that  (\ref{cond-desejavel2}) holds. On the other hand, using (\ref{asstriple}) we obtain
\begin{eqnarray*} \sum_{j=1}^3\delta_jv_j'V_j'&=&\sum_{j=1}^3\delta_jv_j'V_j+\frac{\epsilon\beta}{\varphi}\sum_{j=1}^3\delta_jv_j'v_j-
\frac{\epsilon\beta}{\varphi}\sum_{j=1}^3\delta_j{v_j'}^2\\
&=&\sum_{j=1}^3\delta_jv_j'V_j - \frac{\epsilon\beta}{\varphi}(K -\sum_{j=1}^3\delta_jv_j'v_j)\\
&=&\varphi^{-1}\Omega=0
\end{eqnarray*}
by (\ref{vi'}) and (\ref{omega3}). Thus, the  pair $(v', V')$ associated to  $f'$  also satisfies (\ref{sphol}) (respectively, (\ref{eq:cflatcond})). \qed

\subsection{Explicit three-dimensional solutions of Problem $*$}

We now use Theorem \ref{thm:ribsol} to compute explicit examples of pairs of isometric immersions   $f \colon M^3 \to \Q_s^4(c)$ 
and $\tilde f\colon M^3\to \Q_{\tilde s}^4(\tilde c)$, $c\neq \tilde c$,  with three distinct principal curvatures.

 First notice that, if   $c=0$ (respectively, $c\neq 0$) and $(v,h,V)$ is a solution of system (\ref{sistema-hol}) on a simply connected open subset   $U \subset \R^{3}$ with  $v_i \neq 0$
everywhere  for  $1\leq i\leq 3$, then, in order to determine the corresponding  immersion $f\colon U\to \R_s^4$ (respectively, $f \colon U \to \Q_s^4(c)\subset \R^5_{s+\epsilon_0}$, where $\epsilon_0=c/|c|$), one has to integrate the system of PDE's
\be\label{spde}
 \left\{\begin{array}{l}
  (i) \dfrac{\partial f}{\partial u_i}=v_iX_i,\,\,\,\,\,\,\,
  (ii) \dfrac{\partial X_{i}}{\partial u_j}=h_{ij}X_j,\,\,\,i\neq j,\vspace{1ex}\\
  (iii) \dfrac{\partial X_{i}}{\partial u_i}=-\sum_{k\neq i}h_{ki}X_k+\epsilon V_iN-cv_if,\vspace{1ex}\\
  (iv) \dfrac{\partial N}{\partial u_i}=-V_iX_i, \,\,\,1\leq i\leq 3,
    \end{array}\right.
 \ee
 with initial conditions $X_1(u_0), X_2(u_0), X_3(u_0), N(u_0), f(u_0)$ at some point $u_0\in U$ chosen so that the set $\{X_1(u_0), X_2(u_0), X_3(u_0), N(u_0)\}$  (respectively, 
 $\{X_1(u_0), X_2(u_0), X_3(u_0), N(u_0), |c|^{1/2}f(u_0)\}$) is an orthonormal basis of $\R_s^4$ (respectively, $\R^5_{s+\epsilon_0}$).

   The idea for the construction of explicit examples is to  start with trivial solutions $(v,h,V)$  of system (\ref{sistema-hol}). 
  If $\hat\epsilon=1$, one can start with the solution $(v,h,V)$  of system (\ref{sistema-hol}), with $(\delta_1, \delta_2, \delta_3)=(1, -1, 1)$, for which 
 $v=(1, 0, 0)$, $h=0$ and $V$ is either $\sqrt{-C}(0,1,0)$ or $\sqrt{C}(0,0,1)$, corresponding to $C<0$ or $C>0$, respectively. 
 If $\hat\epsilon=-1$  and $C>0$,  we may start with the  solution $(v,h,V)$  of system (\ref{sistema-hol}),  with $(\delta_1, \delta_2, \delta_3)=(1, -1, 1)$,  for which 
 $v=(0, 1, 0)$, $h=0$ and $V=\sqrt{C}(0,0,1)$, whereas for  $C<0$ we take  $(\delta_1, \delta_2, \delta_3)=(-1, -1, -1)$,  
 $v=(0, 0, 1)$, $h=0$ and $V=\sqrt{-C}(1,0,0)$. Even though, for the corresponding solution $(X_1, X_2, X_3, N, f)$ of system 
 (\ref{spde}), the map $f\colon U\to \Q_s^4(c)$ is not an immersion, the map  $f'\colon U\to \Q_s^4(c)$ obtained by applying
 Theorem \ref{thm:ribsol} to it does define a hypersurface of $\Q_s^4(c)$, which is therefore a solution of Problem~$*$.
 
 In the following, we consider the case in which $\hat\epsilon=1$ and $C<0$, the others being similar. We take 
 $(v,h,V)$ as the solution  of system (\ref{sistema-hol}), with $(\delta_1, \delta_2, \delta_3)=(1, -1, 1)$, for which 
 $v=(1, 0, 0)$, $h=0$ and $V=\sqrt{-C}(0,1,0)$.


 If $c=0$, the corresponding solution of system (\ref{spde}) with initial conditions 
 $$(X_1(0), X_2(0), X_3(0), N(0), f(0))=(E_1,E_2, E_3,\epsilon E_4,0)$$  
 is given by
 $$f=f(u_1)=u_1E_1,\,\,\,\,\,X_1=E_1,\,\,\,\,X_3=E_3,$$
 \be\label{eq:x2b}X_2=\left\{\begin{array}{l}
  \cosh au_2 E_2+\sinh au_2 E_4, \,\,\mbox{if}\,\,\epsilon=-1, \\
  \cos au_2 E_2+\sin au_2 E_4, \,\,\mbox{if}\,\,\epsilon=1, \\
       \end{array}\right.
       \ee
       and 
  \be\label{eq:nb}
  N=\left\{\begin{array}{l}
  -\sinh au_2 E_2-\cosh au_2 E_4, \,\,\mbox{if}\,\,\epsilon=-1, \\
  -\sin au_2 E_2+\cos au_2 E_4, \,\,\mbox{if}\,\,\epsilon=1, \\
       \end{array}\right.
       \ee 
where $a=\sqrt{-C}$.    If $c\neq 0$, the corresponding solution of system (\ref{spde}) with initial conditions 
    $$(X_1(0), X_2(0), X_3(0), N(0), f(0))=(E_1,E_2, E_3,E_4,|c|^{-1/2}E_5)$$
     is given by
 \be\label{fb}
 f=f(u_1)=\left\{\begin{array}{l}
 \frac{1}{\sqrt{c}}(\cos \sqrt{c}\,u_1E_5+\sin \sqrt{c}\,u_1E_1),\,\,\mbox{if}\,\,c>0,\\
 \frac{1}{\sqrt{-c}}(\cosh \sqrt{-c}\,u_1E_5+\sinh \sqrt{-c}\,u_1E_1),\,\,\mbox{if}\,\,c<0,
 \end{array}\right.
       \ee
 \be\label{x1b}
 X_1=\left\{\begin{array}{l}
 -\sin \sqrt{c}\,u_1 E_5 +\cos \sqrt{c}\,u_1E_1,\,\,\mbox{if}\,\,c>0,\\
 \sinh \sqrt{-c}\,u_1 E_5 +\cosh \sqrt{-c}\,u_1E_1,\,\,\mbox{if}\,\,c<0,
 \end{array}\right.
       \ee
 $X_3=E_3$ and  $X_2$, $N$ as in (\ref{eq:x2b}) and (\ref{eq:nb}), respectively.

We now solve  system (\ref{sist-RibA}) for $(v,h,V)$ as in the preceding paragraph. Notice that
(\ref{int-primeiraA}) and (\ref{int-primeiraB}), with $K_1=0$ and $K_2=1$, 
reduce, respectively,  to
\begin{equation}\label{1int1-cfb}
2\varphi \psi = \sum_i\gamma_i^2 + \epsilon\beta^2 +c\varphi^2
\end{equation}
and 
\begin{equation}\label{2int1-cfb}
{v'}_2^2={v'}_1^2+{v'}_3^2-1.
\end{equation}
We also impose that
\begin{eqnarray}\label{impo-cfb}
-a\varphi v_2' = \epsilon\beta(1-v_1'),
\end{eqnarray}
which corresponds to the function $\Omega$ in (\ref{eq:omega}) vanishing everywhere.
It follows from equations $(i)$, $(ii)$ and $(iv)$ of (\ref{sist-RibA}) that $\varphi$, $\gamma_j$ and $\beta$ depend only on $u_1$, $u_j$ and $u_2$, respectively.
Equation $(iii)$ then implies that there exist smooth functions $\phi_i=\phi_i(u_i)$, $1\leq i\leq 3$, such that
\begin{equation}\label{v1'cfb}
(\delta_{1i}-v_i')\psi=\phi_i.
\ee
Replacing (\ref{v1'cfb})  in (\ref{2int1-cfb}) gives
\be\label{eq:psib}
\psi=\frac{\phi_1^2-\phi_2^2+\phi_3^2}{2\phi_1}.
\ee
Multiplying (\ref{impo-cfb}) by $\psi$ and using (\ref{v1'cfb})  yields
$$a\varphi \phi_2=\epsilon \beta \phi_1,$$
hence there exists $K\neq 0$ such that
\be\label{betavar}
\beta=\frac{\epsilon}{K}\phi_2\,\,\,\,\mbox{and}\,\,\,\,\varphi=\frac{1}{Ka}\phi_1.
\ee
It follows from $(i)$ and $(iv)$ that
\be\label{gammai}
\gamma_1=\frac{1}{Ka}\phi_1'\,\,\,\,\mbox{and}\,\,\,\,\gamma_2=-\frac{1}{Ka}\phi_2'
\ee
where $\phi_i'$ stands for the derivative of $\phi_i$ (with respect to $u_i$). Using $(v)$ for $i=3$, (\ref{v1'cfb}) and the second equation in (\ref{betavar}) we obtain that
$$\gamma_3=\frac{1}{Ka}\phi_3'.$$
Then, it follows  from  $(iii)$, (\ref{v1'cfb}), the first equation in (\ref{gammai}) and the second one in (\ref{betavar}) that
\be\label{edophi1}\phi_1''=(Ka-c)\phi_1.\ee
Similarly,
\be\label{edophi23}\phi_2''=-(\epsilon a^2+Ka)\phi_2\,\,\,\,\mbox{and}\,\,\,\,\phi_3''=Ka\phi_3.\ee
Moreover, by (\ref{1int1-cfb}) we must have
\be\label{sumfi}
\phi_1'^2-(Ka-c)\phi_1^2 + \phi_2'^2+(\epsilon a^2+Ka)\phi_2^2 + \phi_3'^2-Ka\phi_3^2=0.
\ee
Notice that each of the expressions under brackets in the preceding equation 
is constant, as follows from (\ref{edophi1}) and (\ref{edophi23}). 

We compute explicitly the corresponding  hypersurface given by (\ref{dt2:10A}) 
when $c=0$, $\tilde c=1$, $\epsilon=1=\tilde \epsilon$ and $K=1$.  In this case we have
$C=-1$ and $a=1$, hence equations (\ref{edophi1}) and (\ref{edophi23}) yield
$$
\left\{\begin{array}{l}
\ \phi_1=A_{11}\cosh u_1+A_{12}\sinh u_1,\\
\ \phi_2=A_{21}\cos \sqrt{2}\,u_2+A_{22}\sin \sqrt{2}\,u_2,\\
\ \phi_3=A_{31}\cosh u_3+A_{32}\sinh u_3,
\end{array}\right.
$$
where $A_{ij}\in \R$, $1\leq i, j\leq 3$, satisfy 
$$A_{12}^2-A_{11}^2+2(A_{21}^2+A_{22}^2)+A_{32}^2-A_{31}^2=0,$$
in view of (\ref{sumfi}). Assuming, say, that 
$$A_{12}^2-A_{11}^2<0\,\,\,\,\mbox{and}\,\,\,\,A_{32}^2-A_{31}^2<0,$$
 we may write $A_{11}=\rho_1\cosh \theta_1$, $A_{12}=\rho_1\sinh \theta_1$,  $A_{21}=\rho_2\sin \theta_2$, $A_{22}=\rho_2\cos  \theta_2$, $A_{31}=\rho_3\cosh \theta_3$ and $A_{32}=\rho_3\sinh \theta_3$ for some $\rho_i>0$ and $\theta_i\in \R$, $1\leq i\leq 3$.  Then 
$$
\left\{\begin{array}{l}
\ \phi_1=\rho_{1}\cosh (u_1+\theta_1),\\
\ \phi_2=\rho_{2}\sin (\sqrt{2}\,u_2+\theta_{2}),\\
\ \phi_3=\rho_{3}\cosh (u_3+\theta_3),
\end{array}\right.
$$
with 
$$
2\rho_2^2=\rho_1^2+\rho_3^2,
$$
and we can assume that $\theta_i=0$ after a suitable change $u_i\mapsto u_i+u_i^0$ of the coordinates $u_i$,  $1\leq i\leq 3$.
Setting $\rho=\rho_2$, we can write $\rho_1=\sqrt{2}\rho\cos\theta$ and $\rho_3=\sqrt{2}\rho\sin\theta$ for some $\theta\in [0, 2\pi]$. Thus
$$
\left\{\begin{array}{l}
\ \phi_1=\sqrt{2}\rho\cos\theta\cosh u_1,\\
\ \phi_2=\rho\sin \sqrt{2}\,u_2,\\
\ \phi_3=\sqrt{2}\rho\sin\theta\cosh u_3,
\end{array}\right.
$$
and the coordinate functions of the corresponding  hypersurface $f'\colon U\to  \R^4$ are
$$f'_1=u_1 - 2 gh\cos\theta\sinh u_1,\,\,\,\,f'_2=gh(2\cos\sqrt{2}u_2\cos u_2+\sqrt{2}\sin\sqrt{2}u_2\sin u_2),$$
$$f'_3=-2gh\sin \theta \sinh u_3,\,\,\,\,\,\, f'_4=gh(2\cos \sqrt{2}u_2\sin u_2- \sqrt{2}\sin \sqrt{2}u_2\phi_2\cos u_2),$$
where
$$g=2\cos \theta\cosh u_1$$
and 
$$h^{-1}=2\cos^2\theta\cosh^2u_1 - \sin^2 \sqrt{2} u_2 + 2\sin^2\theta\cosh^2u_3.$$

To determine the immersion $\tilde{f'}\colon U\to \Sf^4$ that has the same induced metric as $f'$, we start with 
the solution $(\tilde v,\tilde h, \tilde V)$ of  system (\ref{sistema-hol}), together with equations (\ref{eq:int1}) and (\ref{eq:int2}), with $(\delta_1, \delta_2, \delta_3)=(1, -1, 1)$ and $c$ replaced by $\tilde c=1$, for which  $\tilde v=v=(1, 0, 0)$, $\tilde h=h=0$ and $\tilde{V}=(0,0,1)$.

The corresponding solution $(\tilde X_1, \tilde X_2, \tilde X_3, \tilde N, \tilde f)$ of system (\ref{spde}),  with $\epsilon=\tilde \epsilon=1$, $c=\tilde c=1$ and initial conditions 
    $$(\tilde X_1(0), \tilde X_2(0), \tilde X_3(0), \tilde N(0), \tilde f(0))=(E_1,E_2, E_3,E_4,E_5)$$
     is given by
 \be\label{fb1}
 \tilde f=\tilde f(u_1)=\cos u_1E_5+\sin u_1E_1,\ee
 \be\label{x1b1}
 \tilde X_1=-\sin u_1 E_5 +\cos u_1E_1,\;\;\;\;\tilde X_2=E_2,\ee
\be\label{eq:x2c}\tilde X_3=
  \cos u_3 E_3+\sin u_3 E_4\;\;\;\mbox{and}\;\;\;
        \tilde N= -\sin u_3 E_3+\cos u_3 E_4.
       \ee  
 
 Arguing as before, we solve system (\ref{sist-RibA}) together with equations       
 (\ref{int-primeiraA}) and (\ref{int-primeiraB}), which now become
\begin{equation}\label{1int1-cfb1}
2\tilde\varphi \tilde\psi = \sum_i\tilde\gamma_i^2 + \tilde\beta^2 +\tilde\varphi^2
\end{equation}
and 
\begin{equation}\label{2int1-cfb1}
\tilde{v'}_2^2=\tilde{v'}_1^2+\tilde{v'}_3^2-1.
\end{equation}
We also impose that
\begin{eqnarray}\label{impo-cfb1}
\tilde\varphi \tilde v_3' = \tilde\beta(1-\tilde v_1'),
\end{eqnarray}
which corresponds to the function $\Omega$ in (\ref{eq:omega}) vanishing everywhere.
We obtain 
\be\label{eq:psib1}
\tilde\psi=\frac{\tilde\phi_1^2-\tilde\phi_2^2+\tilde\phi_3^2}{2\tilde\phi_1},\;\;\;\;(\delta_{i1}-\tilde v_i')\tilde \psi=\tilde\phi_i,
\ee
\be\label{betavar1}
\tilde \beta=\frac{1}{\tilde K}\tilde\phi_3,\,\,\,\,\,\,\,\tilde\varphi=-\frac{1}{\tilde K}\tilde\phi_1,
\ee
\be\label{gammai1}
\tilde \gamma_1=-\frac{1}{\tilde K}\tilde\phi_1',\,\,\,\,\tilde \gamma_2=\frac{1}{\tilde K}\tilde\phi_2'\,\,\,\mbox{and}\,\,\,\,\tilde \gamma_3=-\frac{1}{\tilde K}\tilde\phi_3'
\ee
for some $\tilde K\in \R$, where the functions $\tilde\phi_i=\tilde\phi_i(u_i)$ satisfy
\be\label{edophi1b}\tilde\phi_1''=-(1+\tilde K)\tilde\phi_1,\,\,\,\,\tilde\phi_2''=\tilde K\tilde\phi_2\,\,\,\,\,\,\tilde\phi_3''=-(1 + \tilde K)\tilde\phi_3\ee
and
\be\label{sumfi1}
(\tilde\phi_1'^2+(1+\tilde K)\tilde\phi_1^2) + (\phi_2'^2-\tilde K\tilde\phi_2^2) + (\phi_3'^2+(1 + \tilde K)\tilde\phi_3^2)=0.
\ee
Notice that each of the expressions under brackets in the preceding equation 
is constant, as follows from (\ref{edophi1b}). Notice also that, for $\tilde K=-2$,  the two preceding equations coincide with (\ref{edophi1}), (\ref{edophi23}) and (\ref{sumfi}) for $1=K=a=\epsilon$ and $c=0$, hence the metrics induced by $f'$ and $\tilde f'\colon U\to \Sf^4\subset \R^5$ coincide by  (\ref{v1'cfb}) and the second equation in (\ref{eq:psib1}). The coordinate functions of $\tilde f'$ are
\begin{eqnarray}
\begin{array}{lcl}
\tilde{f}'_1 &=& \sin u_1 + gh(\cos\theta\cos u_1\sinh u_1+\cos\theta \sin u_1 \cosh u_1)\\
\tilde{f}'_2 &=& -gh\cos \sqrt{2}u_2\\
\tilde{f}'_3 &=& gh(\sin \theta\cos u_3\sinh u_3 - \sin\theta \sin u_3 \cosh u_3)\\
\tilde{f}'_4 &=& gh(\sin \theta\sin u_3\sinh u_3 + \sin\theta \cos u_3 \cosh u_3)\\
\tilde{f}'_5 &=& \cos u_1 + gh(\cos\theta \cos u_1 \cosh u_1 - \cos\theta\sin u_1\sinh u_1)
\end{array}
\end{eqnarray}
where
$$g=2\cos \theta \cosh u_1$$
and 
$$h^{-1}=2\cos^2\theta\cos^2u_1- \sin^2\sqrt{2} u_2 + 2\sin^2\theta\cosh^2u_3.$$

\subsection{Examples of  conformally flat hypersurfaces}

 One can also use Theorem \ref{thm:ribsol} to compute explicit examples of conformally flat hypersurfaces   $f \colon M^3 \to \Q_s^4(c)$ 
with three distinct principal curvatures. It suffices to consider the case $c=0$, because any 
conformally flat hypersurface   $f \colon M^3 \to \Q_s^4(c)$, $c\neq 0$, is the composition of a conformally flat hypersurface
   $g \colon M^3 \to \R_s^4$ with an ``inverse stereographic projection".
   
   We start with the trivial solution $v=(0,1,1)$, $V=(1,0,0)$ and $h=0$ of  system (\ref{sistema-hol}), for which the corresponding
   solution  of system (\ref{spde}) with initial conditions 
 $$(X_1(0), X_2(0), X_3(0), N(0), f(0))=(E_1,E_2, E_3,E_4,0)$$  
 is given by
 $$f=f(u_2, u_3)=u_2E_2+u_3E_3,\,\,\,\,\,X_2=E_2,\,\,\,\,X_3=E_3,$$
 \be\label{eq:x2}X_1=\left\{\begin{array}{l}
  \cosh u_1 E_1+\sinh u_1 E_4, \,\,\mbox{if}\,\,\epsilon=-1,\\
  \cos u_1 E_1+\sin u_1 E_4, \,\,\mbox{if}\,\,\epsilon=1,  
       \end{array}\right.
       \ee
       and 
  \be\label{eq:n}
  N=\left\{\begin{array}{l}
    \sinh u_1 E_1+\cosh u_1 E_4, \,\,\mbox{if}\,\,\epsilon=-1, \\
   -\sin u_1 E_1+\cos u_1 E_4, \,\,\mbox{if}\,\,\epsilon=1. 
       \end{array}\right.
       \ee  
       
Even though this solution does not correspond to a three-dimensional hypersurface, one  can still
  apply Theorem \ref{thm:ribsol}.   We  solve  system (\ref{sist-RibA}) for $(v,h,V)$ as in the preceding paragraph. Equations
(\ref{int-primeiraA}) and (\ref{int-primeiraB}), with $K_1=0=K_2$,   become       
\begin{equation}\label{1int1-cf}
2\varphi \psi = \sum_i\gamma_i^2 + \epsilon\beta^2
\end{equation}
and 
\begin{equation}\label{2int1-cf}
{v'}_2^2={v'}_1^2+{v'}_3^2.
\end{equation}
We also impose that
\begin{eqnarray}\label{impo-cf}
\varphi v_1' = - \epsilon\beta \left(v_3' - v_2'\right),
\end{eqnarray}
which corresponds to the function $\Omega$ in (\ref{eq:omega}) vanishing everywhere.
It follows from $(iii)$ that
$$v_1'\psi=\beta-\frac{\d \gamma_1}{\d u_1}.$$
Since the right-hand-side of the preceding equation  depends only on $u_1$ by $(ii)$ and  $(iv)$, there exists a smooth function $\phi_1=\phi_1(u_1)$ such that 
\begin{equation}\label{v1'cf}
v_1'\psi = \phi_1.
\end{equation}
 Similarly, 
\be\label{v2'v3'}
(1-v_i')\psi=\phi_i
\ee
for some smooth functions $\phi_i=\phi_i(u_i)$, $2\leq i\leq 3$. In particular, 
\be\label{v2v3}
(v_2'-v_3')\psi=\phi_3-\phi_2.
\ee
Multiplying (\ref{impo-cf}) by $\psi$ and using (\ref{v1'cf}) and (\ref{v2v3}) yields
\be\label{varphi}
\varphi=\frac{1}{K}(\phi_3-\phi_2)
\ee
and 
\be\label{beta}
\beta=\frac{\epsilon}{K}\phi_1
\ee
for some $K\in \R$. On the other hand, replacing (\ref{v1'cf}) and (\ref{v2'v3'}) in (\ref{2int1-cf}) gives
\be\label{eq:psi}
\psi=\frac{\phi_1^2-\phi_2^2+\phi_3^2}{2(\phi_3-\phi_2)}.
\ee
It follows from $(i)$ and $(\ref{varphi})$ that  
$$\gamma_2=-\frac{1}{K}\phi_2'\,\,\,\,\mbox{and}\,\,\,\,\,\gamma_3=\frac{1}{K}\phi_3',$$
whereas $(iv)$ and (\ref{beta}) yields
\be\label{gamma1}
\gamma_1=-\frac{1}{K}\phi_1'.
\ee
We obtain from $(iii)$, (\ref{beta}) and (\ref{gamma1}) that
\be\label{f1}
\phi_1''=(K-\epsilon)\phi_1.
\ee
Similarly, 
\be\label{f2f3}
\phi_2''=-K\phi_2\,\,\,\,\mbox{and}\,\,\,\,\phi_3''=K\phi_3.
\ee
Moreover, by (\ref{1int1-cf}) we must have
\be\label{sumfib}
(\phi_1'^2-(K-\epsilon)\phi_1^2) + (\phi_2'^2+K\phi_2^2) + (\phi_3'^2-K\phi_3^2)=0.
\ee
Notice that each of the expressions under brackets in the preceding equation 
is constant, as follows from (\ref{f1}) and (\ref{f2f3}). 

The conformally flat hypersurface given by (\ref{dt2:10A}) (with $c=0$) has coordinate functions
$$ f'_1=-\psi^{-1}(\phi_1'\cos u_1+\phi_1\sin u_1),\,\,\,\, f'_2=u_2-\psi^{-1} \phi_2',$$
$$ f'_3=u_3+\psi^{-1} \phi_3'\,\,\,\,\mbox{and}\,\,\,\, f'_4=\psi(\phi_1\cos u_1-\phi_1'\sin u_1),$$
with $\psi$  as in (\ref{eq:psi}).  We compute them explicitly  for the particular case
$\epsilon=1$ and $K<0$, the others being similar. In this case we have
$$
\left\{\begin{array}{l}
\ \phi_1=A_{11}\cos \sqrt{|K-1|}\,u_1+A_{12}\sin \sqrt{|K-1|}\,u_1,\\
\ \phi_2=A_{21}\cosh \sqrt{|K|}\,u_2+A_{22}\sinh \sqrt{|K|}\,u_2,\\
\ \phi_3=A_{31}\cos \sqrt{|K|}\,u_3+A_{32}\sin \sqrt{|K|}\,u_3,
\end{array}\right.
$$
with $A_{ij}\in \R$ for $1\leq i, j\leq 3$, and equation (\ref{sumfib}) reduces to
$$|K-1|(A_{11}^2+A_{12}^2)+|K|(A_{22}^2-A_{21}^2)+|K|(A_{31}^2+A_{32}^2)=0.$$
This implies that
$$A_{22}^2-A_{21}^2<0,$$
hence we may write $A_{21}=\rho_2\cosh \theta_2$ and $A_{22}=\rho_2\sinh \theta_2$ for some $\rho_2>0$ and $\theta_2\in \R$. 
We may also write $A_{11}=\rho_1\cos \theta_1$, $A_{12}=\rho_1\sin \theta_1$, $A_{31}=\rho_3\cos \theta_3$ and $A_{32}=\rho_3\sin \theta_3$ for some $\rho_1, \rho_3>0$ and $\theta_1, \theta_3\in [0, 2\pi]$. Then 
$$
\left\{\begin{array}{l}
\ \phi_1=\rho_{1}\cos (\sqrt{|K-1|}\,u_1-\theta_1),\\
\ \phi_2=\rho_{2}\cosh (\sqrt{|K|}\,u_2+\theta_{2}),\\
\ \phi_3=\rho_{3}\cos (\sqrt{|K|}\,u_3-\theta_3),
\end{array}\right.
$$
with 
$$
|K|\rho_2^2=|K-1|\rho_1^2+|K|\rho_3^2,
$$
and we can assume that $\theta_i=0$ after a suitable change $u_i\mapsto u_i+u_i^0$ of the coordinates $u_i$,  $1\leq i\leq 3$.
Setting $\rho=\rho_2$, we can write $\rho_1=\sqrt{\frac{|K|}{|K-1|}}\rho\cos\theta$ and $\rho_3=\rho\sin\theta$ for some $\theta\in [0, 2\pi]$. Thus
$$
\left\{\begin{array}{l}
\ \phi_1=\sqrt{\frac{|K|}{|K-1|}}\rho\cos\theta\cos (\sqrt{|K-1|}\,u_1),\\
\ \phi_2=\rho\cosh (\sqrt{|K|}\,u_2),\\
\ \phi_3=\rho\sin\theta\cos (\sqrt{|K|}\,u_3).
\end{array}\right.
$$
For instance, for $K=-1$ we get the conformally flat hypersurface of $\R^4$ whose coordinate functions are
$$f'_1=2\cos\theta g h (\sqrt{2}\cos \sqrt{2}u_1 \sin u_1-\sin \sqrt{2}u_1 \cos u_1), \,\,\,\,f'_2=u_2+4\sinh u_2 g h,$$ 
$$f'_3=u_3+4\sin\theta\sin u_3 g h, \,\,\,\,f'_4=-2\cos \theta(\sin \sqrt{2}u_1\sin u_1+\cos\sqrt{2}u_1\sin u_1)g h$$
where
$$g=\cosh u_2-\sin \theta\cos u_3$$
and 
$$h^{-1}=\cos^2\theta\cos^2\sqrt{2}u_1-2\cosh^2u_2 + 2\sin^2\theta\cos^2u_3.$$

\subsection{Proof of Corollary \ref{parallel}}

Given a hypersurface  $f\colon M^3 \to \Q_s^4(c)$, set $\epsilon =-2s+1$, $\epsilon_c=c/|c|$ and $\check \epsilon=\epsilon\epsilon_c$. 
Let $\varphi$ and $\psi$ be defined by 
$$(\varphi(t), \psi(t))=\left\{\begin{array}{l}(\cos (\sqrt{|c|} t), \sin(\sqrt{|c|} t)),\,\,\,\mbox{if}\,\,\,\check \epsilon=1,\vspace{1ex}\\
(\cosh (\sqrt{|c|} t), \sinh(\sqrt{|c|} t)), \,\,\,\mbox{if} \,\,\,\check \epsilon=-1.
\end{array}\right.
$$
Then the family of parallel hypersurfaces $f_t\colon M^3 \to \Q_s^4(c)\subset \R_{s+\epsilon_0}^5$ to $f$ is given by 
$$i\circ f_t=\varphi(t)i\circ f+\frac{\psi(t)}{\sqrt{|c|}}i_*N,$$
where $N$ is one of the unit normal vector fields to $f$ and $i\colon \Q_s^4(c)\to \R_{s+\epsilon_0}^5$ is the inclusion, with $\epsilon_0=0$ or $1$, corresponding to $c>0$ or $c<0$, respectively.   We denote by $M_t^3$ the manifold $M^3$ endowed with the metric induced by $f_t$.

\begin{proposition}\label{parallel2}
Let  $f\colon M^3 \to \Q_s^4(c)$ be a holonomic hypersurface. Then any parallel hypersurface
$f_t\colon M_t^3 \to \Q_s^4(c)$ to $f$ is also  holonomic and
the  pairs $(v,  V)$ and $(v^t,  V^t)$ associated to $f$ and $f_t$, respectively, are related by 
\be\label{vit}\left\{\begin{array}{l}v^t_i=\varphi(t) v_i-\frac{\psi(t)}{\sqrt{|c|}}V_i\vspace{1ex}\\
V_i^t=\check\epsilon \sqrt{|c|}\psi(t) v_i+\varphi(t) V_i.\end{array}\right.
\ee
In particular, $h_{ij}^t=h_{ij}$.
\end{proposition} 
\proof We have
\be\label{ft}
{f_t}_*=\varphi(t)f_*+\frac{\psi(t)}{\sqrt{|c|}}N_*=f_*\left(\varphi(t)I-\frac{\psi(t)}{\sqrt{|c|}}A\right),
\ee
thus a unit normal vector field to $f_t$ is
$$N_t=-\check \epsilon\sqrt{|c|}\psi(t)f+\varphi(t)N.$$
Then,
\begin{eqnarray*}{N_t}_*&=&f_*\left(-\check \epsilon\sqrt{|c|}\psi(t) I-\varphi(t)A\right)\\
&=&-{f_t}_*\left(\varphi(t)I-\frac{\psi(t)}{\sqrt{|c|}}A\right)^{-1}\left(\check \epsilon\sqrt{|c|}\psi(t) I+\varphi(t)A\right).
\end{eqnarray*}
which implies that
\be\label{at}
A_t=\left(\varphi(t)I-\frac{\psi(t)}{\sqrt{|c|}}A\right)^{-1}\left(\check \epsilon\sqrt{|c|}\psi(t) I+\varphi(t)A\right).
\ee
It follows from (\ref{ft}) and (\ref{at}) that $\tilde f$ is also holonomic with associated pair given by (\ref{vit}). The assertion on  $h_{ij}^t$ follows 
from a straightforward computation.\vspace{2ex}\qed\\
\noindent \emph{Proof of Corollary \ref{parallel}:} Conditions (\ref{sphol}) for $(v^t,  V^t)$ (with $\tilde c=0$) follow immediately from those for $(v,V)$. \qed
\begin{remark}\emph{Given a hypersurface  $f\colon M^3 \to \Q_s^4(c)$, it can be checked that the 
parallel hypersurfaces $f_t\colon M^3 \to \Q_s^4(c)$ correspond to the Ribaucour transforms of $f$ determined by solutions
 $(\gamma_1,\gamma_2,\gamma_3,v_1',v_2',v_3',\varphi,\psi,\beta)$ of system (\ref{sist-Rib}) for which $\gamma_1=\gamma_2=\gamma_3=0$ and $\varphi, \psi$ and $\beta$ are constants satisfying (\ref{int-primeirab}).}\end{remark} 

\section{Proof of Theorem \ref{intersection}}

For the proof of Theorem \ref{intersection} we need the following preliminary fact, which was already observed in \cite{dt1} for $s=0$.

\begin{lemma}\label{le:threecond} Let $f\colon M^3 \to \Q_s^4(c)$ be a hypersurface with three distinct principal curvatures $\lambda_1, \lambda_2$ and $\lambda_3$. Then, any two of the following three conditions imply the remaining one:
\begin{itemize}
\item[$(i)$] $(\lambda_j - \lambda_k)e_i(\lambda_i)+(\lambda_i - \lambda_k)e_i(\lambda_j)+(\lambda_j - \lambda_i)e_i(\lambda_k)=0$. 
\item[$(ii)$] 
$(C+\hat \epsilon\lambda_j\lambda_k)(\lambda_k-\lambda_j)e_i(\lambda_i) + (C+\hat \epsilon\lambda_i\lambda_k)(\lambda_k-\lambda_i)e_i(\lambda_j)\vspace{1ex}\\ \hspace*{20ex}+ (C+\hat \epsilon\lambda_i\lambda_j)(\lambda_i-\lambda_j)e_i(\lambda_k)=0. 
$
\item[$(iii)$] $e_i(\lambda_i\lambda_j)=0$.
\end{itemize}
\end{lemma}
\proof It is easily checked that $(i)$ is equivalent to 
$$(\lambda_k-\lambda_i)e_i(\lambda_i\lambda_j)=(\lambda_j-\lambda_i)e_i(\lambda_i\lambda_k),\,\,\,\,1\leq i\neq j\neq k\neq i\leq 3,$$
whereas the difference between $(ii)$ and $(i)$ is equivalent to 
 $$\lambda_k(\lambda_k-\lambda_i)e_i(\lambda_i\lambda_j)=\lambda_j(\lambda_j-\lambda_i)e_i(\lambda_i\lambda_k),\,\,\,\,1\leq i\neq j\neq k\neq i\leq 3,$$
 and the statement follows easily.\vspace{2ex}\qed\\
\noindent \emph{Proof of Theorem \ref{intersection}:} By Theorem \ref{main3}, $f$ is locally a holonomic hypersurface whose associated pair 
$(v,V)$ is given in terms of  the principal curvatures $\lambda_1, \lambda_2$  and $\lambda_3$ of  $f$ by
\be\label{eq:vj}v_j=\sqrt{\frac{\delta_j}{(\lambda_j-\lambda_i)(\lambda_j-\lambda_k)}}\,\,\,\,\mbox{and}\,\,\,\,V_j=\lambda_jv_j,\,\,\,\,\,\,1\leq j\leq 3.
\ee
Moreover,  we have seen in the proofs of Theorems \ref{main2} and \ref{main3} that conditions $(i)$ and $(ii)$ in Lemma \ref{le:threecond} hold for $\lambda_1, \lambda_2$  and $\lambda_3$. Thus, also condition $(iii)$ is satisfied. Assuming that $\lambda_j\neq 0$ for $1\leq j\leq 3$, we can write
\be\label{lilj}\lambda_i\lambda_j=\iota_k\phi_k^2, \,\,\,\,\, \iota_k\in \{-1, 1\},\,\,\,\,1\leq i\neq j\neq k\neq i\leq 3,\ee
for some positive smooth functions $\phi_k=\phi_k(u_k)$, $1\leq k\leq 3$. It follows from (\ref{lilj}) that 
\be\label{lambdaj}
\lambda_j=\epsilon_j\frac{\phi_i\phi_k}{\phi_j},
\ee
where $\epsilon_j=\frac{\lambda_j}{|\lambda_j|}$, $1\leq j\leq 3$. We may suppose that $\lambda_1<\lambda_2<\lambda_3$, so that
$$\epsilon_k\phi_i^2-\epsilon_i\phi_k^2>0,\,\,\,\,1\leq i<k\leq 3.$$
Substituting (\ref{lambdaj}) into (\ref{eq:vj}), we obtain that
\be\label{eq:vj2}v_j=\frac{\phi_j}{\psi_i\psi_k},\,\,\,\,\,\,1\leq j\leq 3,
\ee
where 
$$\psi_j=\sqrt{\epsilon_k\phi_i^2-\epsilon_i\phi_k^2}$$
and
$$V_j=\lambda_jv_j=\epsilon_j\frac{\phi_i\phi_k}{\psi_i\psi_k},\,\,\,\,\,\, i, k\neq j, \,\,\,\,\, i<k.$$
We obtain from (\ref{eq:vj2}) that
\be\label{hij1}
h_{ij}=\frac{1}{v_j}\frac{\d v_j}{\d u_i}=\frac{\psi_i\psi_k}{\phi_j}\frac{\phi_j}{\psi_i\psi_k^2}\left(-\frac{\d \psi_k}{\d u_i}\right)=-\frac{1}{\psi_k}\frac{\d \psi_k}{\d u_i}.
\ee
On the other hand, equation $(iv)$ of system (\ref{sistema-hol}) yields
\be\label{hij2}
h_{ij}=\frac{1}{V_j}\frac{\d V_j}{\d u_i}=\frac{\psi_i\psi_k}{\phi_i\phi_k}\frac{\phi_k}{\psi_i\psi_k^2}\left(\frac{d \phi_i}{d u_i}\psi_k-\phi_i\frac{\d \psi_k}{\d u_i}\right)=\frac{1}{\phi_i}\frac{d \phi_i}{d u_i}-\frac{1}{\psi_k}\frac{\d \psi_k}{\d u_i}.
\ee
Comparying (\ref{hij1}) and (\ref{hij2}), we obtain that
$$\frac{d \phi_i}{d u_i}=0,\,\,\,\,1\leq i\leq 3.$$
This implies that $\frac{\d \psi_k}{\d u_i}=0$ for all $1\leq i\neq k\leq 3$, and hence $h_{ij}=0$ for all $1\leq i\neq j\leq 3$. 
But then equation $(ii)$ of system (\ref{sistema-hol}) gives
$$\epsilon\lambda_i\lambda_j+c=0$$
for all $1\leq i\neq j\leq 3$, which implies that $-\epsilon c>0$ and $\lambda_1=\lambda_2=\lambda_3=\sqrt{-\epsilon c}$, a contradiction. 
Thus, one of the principal curvatures must be zero, and the result follows from part $b)$ of Theorem \ref{multdoisb}.\qed

{\renewcommand{\baselinestretch}{1} \hspace*{-20ex}\begin{tabbing}
\indent \= Samuel Canevari\\
\> Universidade Federal de Sergipe \\
\> Av. Vereador Olimpio Grande s/n.\\
\> Itabaiana  -- Brazil\\
\> scanevari@gmail.com\\
\end{tabbing}}

\vspace*{-3ex}

{\renewcommand{\baselinestretch}{1} \hspace*{-20ex}\begin{tabbing}
\indent \= Ruy Tojeiro\\
\> Departamento de Matematica, \\
\> Universidade Federal de S\~{a}o Carlos,\\
\> Via Washington Luiz km 235\\
\> 13565-905 -- S\~{a}o Carlos -- Brazil\\
\> tojeiro@dm.ufscar.br
\end{tabbing}}

\end{document}